\documentclass[reqno,12pt,a4paper]{amsart}

\usepackage{mathrsfs}
\usepackage{amssymb}
\usepackage{amsfonts}
\usepackage{latexsym}
\usepackage{amsthm}
\usepackage{graphicx}
\def\lb{\label}

\newcommand{\er}[1]{\textrm{(\ref{#1})}}

\begin{document}


\renewcommand{\theequation}{\arabic{section}.\arabic{equation}}
\theoremstyle{plain}
\newtheorem{theorem}{\bf Theorem}[section]
\newtheorem{lemma}[theorem]{\bf Lemma}
\newtheorem{corollary}[theorem]{\bf Corollary}
\newtheorem{proposition}[theorem]{\bf Proposition}
\newtheorem{definition}[theorem]{\bf Definition}
\newtheorem{remark}[theorem]{\bf Remark}

\def\a{\alpha}  \def\cA{{\mathcal A}}     \def\bA{{\bf A}}  \def\mA{{\mathscr A}}
\def\b{\beta}   \def\cB{{\mathcal B}}     \def\bB{{\bf B}}  \def\mB{{\mathscr B}}
\def\g{\gamma}  \def\cC{{\mathcal C}}     \def\bC{{\bf C}}  \def\mC{{\mathscr C}}
\def\G{\Gamma}  \def\cD{{\mathcal D}}     \def\bD{{\bf D}}  \def\mD{{\mathscr D}}
\def\d{\delta}  \def\cE{{\mathcal E}}     \def\bE{{\bf E}}  \def\mE{{\mathscr E}}
\def\D{\Delta}  \def\cF{{\mathcal F}}     \def\bF{{\bf F}}  \def\mF{{\mathscr F}}
\def\c{\chi}    \def\cG{{\mathcal G}}     \def\bG{{\bf G}}  \def\mG{{\mathscr G}}
\def\z{\zeta}   \def\cH{{\mathcal H}}     \def\bH{{\bf H}}  \def\mH{{\mathscr H}}
\def\e{\eta}    \def\cI{{\mathcal I}}     \def\bI{{\bf I}}  \def\mI{{\mathscr I}}
\def\p{\psi}    \def\cJ{{\mathcal J}}     \def\bJ{{\bf J}}  \def\mJ{{\mathscr J}}
\def\vT{\Theta} \def\cK{{\mathcal K}}     \def\bK{{\bf K}}  \def\mK{{\mathscr K}}
\def\k{\kappa}  \def\cL{{\mathcal L}}     \def\bL{{\bf L}}  \def\mL{{\mathscr L}}
\def\l{\lambda} \def\cM{{\mathcal M}}     \def\bM{{\bf M}}  \def\mM{{\mathscr M}}
\def\L{\Lambda} \def\cN{{\mathcal N}}     \def\bN{{\bf N}}  \def\mN{{\mathscr N}}
\def\m{\mu}     \def\cO{{\mathcal O}}     \def\bO{{\bf O}}  \def\mO{{\mathscr O}}
\def\n{\nu}     \def\cP{{\mathcal P}}     \def\bP{{\bf P}}  \def\mP{{\mathscr P}}
\def\r{\rho}    \def\cQ{{\mathcal Q}}     \def\bQ{{\bf Q}}  \def\mQ{{\mathscr Q}}
\def\s{\sigma}  \def\cR{{\mathcal R}}     \def\bR{{\bf R}}  \def\mR{{\mathscr R}}
\def\S{\Sigma}  \def\cS{{\mathcal S}}     \def\bS{{\bf S}}  \def\mS{{\mathscr S}}
\def\t{\tau}    \def\cT{{\mathcal T}}     \def\bT{{\bf T}}  \def\mT{{\mathscr T}}
\def\f{\phi}    \def\cU{{\mathcal U}}     \def\bU{{\bf U}}  \def\mU{{\mathscr U}}
\def\F{\Phi}    \def\cV{{\mathcal V}}     \def\bV{{\bf V}}  \def\mV{{\mathscr V}}
\def\P{\Psi}    \def\cW{{\mathcal W}}     \def\bW{{\bf W}}  \def\mW{{\mathscr W}}
\def\o{\omega}  \def\cX{{\mathcal X}}     \def\bX{{\bf X}}  \def\mX{{\mathscr X}}
\def\x{\xi}     \def\cY{{\mathcal Y}}     \def\bY{{\bf Y}}  \def\mY{{\mathscr Y}}
\def\X{\Xi}     \def\cZ{{\mathcal Z}}     \def\bZ{{\bf Z}}  \def\mZ{{\mathscr Z}}
\def\be{{\bf e}}
\def\bv{{\bf v}} \def\bu{{\bf u}}
\def\Om{\Omega} \def\mn{\mathrm n}
\def\mm{\mathrm m}

\newcommand{\mc}{\mathscr {c}}

\newcommand{\gA}{\mathfrak{A}}          \newcommand{\ga}{\mathfrak{a}}
\newcommand{\gB}{\mathfrak{B}}          \newcommand{\gb}{\mathfrak{b}}
\newcommand{\gC}{\mathfrak{C}}          \newcommand{\gc}{\mathfrak{c}}
\newcommand{\gD}{\mathfrak{D}}          \newcommand{\gd}{\mathfrak{d}}
\newcommand{\gE}{\mathfrak{E}}
\newcommand{\gF}{\mathfrak{F}}           \newcommand{\gf}{\mathfrak{f}}
\newcommand{\gG}{\mathfrak{G}}           
\newcommand{\gH}{\mathfrak{H}}           \newcommand{\gh}{\mathfrak{h}}
\newcommand{\gI}{\mathfrak{I}}           \newcommand{\gi}{\mathfrak{i}}
\newcommand{\gJ}{\mathfrak{J}}           \newcommand{\gj}{\mathfrak{j}}
\newcommand{\gK}{\mathfrak{K}}            \newcommand{\gk}{\mathfrak{k}}
\newcommand{\gL}{\mathfrak{L}}            \newcommand{\gl}{\mathfrak{l}}
\newcommand{\gM}{\mathfrak{M}}            \newcommand{\gm}{\mathfrak{m}}
\newcommand{\gN}{\mathfrak{N}}            \newcommand{\gn}{\mathfrak{n}}
\newcommand{\gO}{\mathfrak{O}}
\newcommand{\gP}{\mathfrak{P}}             \newcommand{\gp}{\mathfrak{p}}
\newcommand{\gQ}{\mathfrak{Q}}             \newcommand{\gq}{\mathfrak{q}}
\newcommand{\gR}{\mathfrak{R}}             \newcommand{\gr}{\mathfrak{r}}
\newcommand{\gS}{\mathfrak{S}}              \newcommand{\gs}{\mathfrak{s}}
\newcommand{\gT}{\mathfrak{T}}             \newcommand{\gt}{\mathfrak{t}}
\newcommand{\gU}{\mathfrak{U}}             \newcommand{\gu}{\mathfrak{u}}
\newcommand{\gV}{\mathfrak{V}}             \newcommand{\gv}{\mathfrak{v}}
\newcommand{\gW}{\mathfrak{W}}             \newcommand{\gw}{\mathfrak{w}}
\newcommand{\gX}{\mathfrak{X}}               \newcommand{\gx}{\mathfrak{x}}
\newcommand{\gY}{\mathfrak{Y}}              \newcommand{\gy}{\mathfrak{y}}
\newcommand{\gZ}{\mathfrak{Z}}             \newcommand{\gz}{\mathfrak{z}}

\def\ve{\varepsilon}   \def\vt{\vartheta}    \def\vp{\varphi}    \def\vk{\varkappa}

\def\A{{\mathbb A}} \def\B{{\mathbb B}} \def\C{{\mathbb C}}
\def\dD{{\mathbb D}} \def\E{{\mathbb E}} \def\dF{{\mathbb F}} \def\dG{{\mathbb G}} \def\H{{\mathbb H}}\def\I{{\mathbb I}} \def\J{{\mathbb J}} \def\K{{\mathbb K}} \def\dL{{\mathbb L}}\def\M{{\mathbb M}} \def\N{{\mathbb N}} \def\O{{\mathbb O}} \def\dP{{\mathbb P}} \def\R{{\mathbb R}}\def\S{{\mathbb S}} \def\T{{\mathbb T}} \def\U{{\mathbb U}} \def\V{{\mathbb V}}\def\W{{\mathbb W}} \def\X{{\mathbb X}} \def\Y{{\mathbb Y}} \def\Z{{\mathbb Z}}


\def\la{\leftarrow}              \def\ra{\rightarrow}            \def\Ra{\Rightarrow}
\def\ua{\uparrow}                \def\da{\downarrow}
\def\lra{\leftrightarrow}        \def\Lra{\Leftrightarrow}


\def\lt{\biggl}                  \def\rt{\biggr}
\def\ol{\overline}               \def\wt{\widetilde}
\def\no{\noindent}


\let\ge\geqslant                 \let\le\leqslant
\def\lan{\langle}                \def\ran{\rangle}
\def\/{\over}                    \def\iy{\infty}
\def\sm{\setminus}               \def\es{\emptyset}
\def\ss{\subset}                 \def\ts{\times}
\def\pa{\partial}                \def\os{\oplus}
\def\om{\ominus}                 \def\ev{\equiv}
\def\iint{\int\!\!\!\int}        \def\iintt{\mathop{\int\!\!\int\!\!\dots\!\!\int}\limits}
\def\el2{\ell^{\,2}}             \def\1{1\!\!1}
\def\sh{\sharp}
\def\wh{\widehat}
\def\bs{\backslash}
\def\intl{\int\limits}

\def\na{\mathop{\mathrm{\nabla}}\nolimits}
\def\sh{\mathop{\mathrm{sh}}\nolimits}
\def\ch{\mathop{\mathrm{ch}}\nolimits}
\def\where{\mathop{\mathrm{where}}\nolimits}
\def\all{\mathop{\mathrm{all}}\nolimits}
\def\as{\mathop{\mathrm{as}}\nolimits}
\def\Area{\mathop{\mathrm{Area}}\nolimits}
\def\arg{\mathop{\mathrm{arg}}\nolimits}
\def\const{\mathop{\mathrm{const}}\nolimits}
\def\det{\mathop{\mathrm{det}}\nolimits}
\def\diag{\mathop{\mathrm{diag}}\nolimits}
\def\diam{\mathop{\mathrm{diam}}\nolimits}
\def\dim{\mathop{\mathrm{dim}}\nolimits}
\def\dist{\mathop{\mathrm{dist}}\nolimits}
\def\Im{\mathop{\mathrm{Im}}\nolimits}
\def\Iso{\mathop{\mathrm{Iso}}\nolimits}
\def\Ker{\mathop{\mathrm{Ker}}\nolimits}
\def\Lip{\mathop{\mathrm{Lip}}\nolimits}
\def\rank{\mathop{\mathrm{rank}}\limits}
\def\Ran{\mathop{\mathrm{Ran}}\nolimits}
\def\Re{\mathop{\mathrm{Re}}\nolimits}
\def\Res{\mathop{\mathrm{Res}}\nolimits}
\def\res{\mathop{\mathrm{res}}\limits}
\def\sign{\mathop{\mathrm{sign}}\nolimits}
\def\span{\mathop{\mathrm{span}}\nolimits}
\def\supp{\mathop{\mathrm{supp}}\nolimits}
\def\Tr{\mathop{\mathrm{Tr}}\nolimits}
\def\BBox{\hspace{1mm}\vrule height6pt width5.5pt depth0pt \hspace{6pt}}


\newcommand\nh[2]{\widehat{#1}\vphantom{#1}^{(#2)}}
\def\dia{\diamond}

\def\Oplus{\bigoplus\nolimits}



\def\qqq{\qquad}
\def\qq{\quad}
\let\ge\geqslant
\let\le\leqslant
\let\geq\geqslant
\let\leq\leqslant
\newcommand{\ca}{\begin{cases}}
\newcommand{\ac}{\end{cases}}
\newcommand{\ma}{\begin{pmatrix}}
\newcommand{\am}{\end{pmatrix}}
\renewcommand{\[}{\begin{equation}}
\renewcommand{\]}{\end{equation}}
\def\eq{\begin{equation}}
\def\qe{\end{equation}}
\def\[{\begin{equation}}
\def\bu{\bullet}
\def\bq{\mathbf q}

\title[{Spectral band bracketing for Laplacians on periodic metric graphs}]
{Spectral band bracketing for Laplacians on periodic metric graphs}

\date{\today}

\author[Evgeny Korotyaev]{Evgeny Korotyaev}
\address{Mathematical Physics Department, Faculty of Physics, Ulianovskaya 2,
St. Petersburg State University, St. Petersburg, 198904,
\ korotyaev@gmail.com,}
\author[Natalia Saburova]{Natalia Saburova}
\address{Department of Mathematical Analysis, Algebra and Geometry, Uritskogo St. 68, Northern (Arctic)
Federal University, Arkhangelsk, 163002, Russia,
 \ n.saburova@gmail.com}


\begin{abstract}
We consider Laplacians on periodic metric graphs with unit-length edges.
The spectrum of these operators consists of an
absolutely continuous part (which is a union of an infinite number
of non-degenerated spectral bands) plus an infinite number of flat
bands, i.e., eigenvalues of infinite multiplicity. Our main result is a localization of spectral bands in terms of eigenvalues of
Dirichlet and Neumann operators on a fundamental domain of the periodic
graph. The proof is based on the spectral band localization for
discrete Laplacians and on the relation between the spectra of discrete and metric Laplacians.
\end{abstract}

\maketitle


\vskip 0.25cm

\section {Introduction and main results}
\setcounter{equation}{0}

We consider Laplace operators on $\Z^d$-periodic metric graphs with unit-length edges, i.e., on the so-called $\Z^d$-periodic equilateral graphs, $d\geq2$. Differential operators on metric graphs
arise naturally as simplified models in mathematics, physics,
chemistry, and engineering.
It is well-known that the spectrum of the Laplacian on periodic
metric graphs consists of an absolutely continuous part plus an
infinite number of flat bands (i.e., eigenvalues with infinite
multiplicity). The absolutely continuous spectrum is a union of an
infinite number of spectral bands separated by gaps.

For the case of periodic metric graphs we
know only two papers about estimates of the bands and gaps:

(1) Lled\'o and Post \cite{LP08} considered the Laplacian on periodic metric graphs. They estimated the position of the spectral bands of the Laplacian in terms of eigenvalues of the Dirichlet and Neumann operators on a fundamental domain of the periodic
graph. Then using the Cattaneo correspondence \cite{C97} between the spectra of discrete and metric Laplacians they carried over this estimate from the metric Laplacian to the discrete one.

(2) Korotyaev and Saburova \cite{KS14a} obtained another type of the estimate for the metric Laplacian on periodic graphs. They estimated the total length of the spectral bands on a finite interval in terms of geometric parameters of the graph only. In order to do this they estimated the Lebesgue measure of the spectrum for the discrete Laplacian on graphs. After this using the Cattaneo correspondence they carried over the estimate from the discrete case to the Laplacian on metric graphs.

Our main goal is to estimate the position of the spectral bands for the Laplacian on equilateral metric graphs using Dirichlet-Neumann bracketing. Our approach is opposite to Lled\'o -- Post's one. They directly estimated the spectral band positions for the metric Laplacian. Then using the Cattaneo correspondence \cite{C97} between the spectra of discrete and metric Laplacians they determined
Dirichlet-Neumann bracketing for the  normalized Laplacian on
periodic  discrete graphs. Finally, they wrote {\it "It is a priori
not clear how the eigenvalue bracketing can be seen directly for
discrete Laplacians, so our analysis may serve as an example of how
to use metric graphs to obtain results for discrete graphs"} (p.809
in \cite{LP08}). Our approach is opposite and is based on the analysis of the discrete Laplacian on graphs. Then we use the Cattaneo correspondence and carry over the spectral band localization for discrete Laplacians to the metric one.

\subsection{Metric Laplacians.} Let $\G=(V,\cE)$ be a connected infinite graph,
possibly  having loops and multiple edges, where $V$ is the set of
its vertices and  $\cE$ is the set of its unoriented edges. The
graphs under consideration are embedded into $\R^d$. An edge
connecting vertices $u$ and $v$ from $V$ will be denoted as the
unordered pair $(u,v)_e\in\cE$ and is said to be \emph{incident} to
the vertices. Vertices $u,v\in V$ will be called \emph{adjacent} and
denoted by $u\sim v$, if $(u,v)_e\in \cE$. We define the degree
${\vk}_v$ of the vertex $v\in V$ as the number of all its
incident edges from $\cE$ (here a loop is counted twice). Below we consider locally finite
$\Z^d$-periodic metric equilateral graphs $\G$, i.e., graphs satisfying the
following conditions:

1) {\it the number of vertices from $V$ in any bounded domain $\ss\R^d$ is
finite;

2) the degree of each vertex is finite;

3) there exists a basis $a_1,\ldots,a_d$ in $\R^d$ such that $\G$
 is invariant under translations through the vectors $a_1,\ldots,a_d$:
$$
\G+a_s=\G, \qqq  \forall\,s\in\N_d=\{1,\ldots,d\}.
$$
The vectors $a_1,\ldots,a_d$ are called the periods of $\G$.}

4) \emph{All edges of the graph have the unit length.}

From this definition it follows that a $\Z^d$-periodic graph $\G$
 is invariant under translations through any integer vector
  (in the basis $a_1,\ldots,a_d$):
$$
\G+\mm=\G,\qqq \forall\, \mm\in\Z^d.
$$

Each edge $\be$ of $\G$ will be identified with the segment $[0,1]$.
This identification introduces a local coordinate $t\in[0,1]$ along
each edge. Thus, we give an orientation on the edge. Note that the
spectrum of Laplacians on metric graphs does not depend on the
orientation of graph edges. For each function $y$ on $\G$ we define
a function $y_{\be}=y\big|_{\be}$, $\be\in\cE$. We identify each
function  $y_{\be}$ on $\be$ with a function on $[0,1]$ by using the
local coordinate $t\in[0,1]$. Let $L^2(\G)$ be the Hilbert space of all functions $y=(y_\be)_{\be\in\cE}$, where each $y\big|_{\be}\in L^2(0,1)$, equipped with the norm
$$
\|y\|^2_{L^2(\G)}=\sum_{\be\in\cE}\|y_\be\|^2_{L^2(0,1)}<\infty.
$$
We define the metric Laplacian $\D_M$ on $y=(y_\be)_{\be\in\cE}\in L^2(\G)$ by
\[\lb{MLa}
(\D_My)_\be=-y''_\be, \qqq (y''_\be)_{\be\in\cE}\in L^2(\G),
\]
where $y$ satisfies the so-called Kirchhoff conditions:
\[
\lb{Dom1}
y \textrm{ is continuous on }\G,\qqq
\sum\limits_{\be=(v,\,u)_e\in\cE}\d_\be(v)\,y_\be'(v)=0, \qq \forall v\in V,
\]
$$
\d_\be(v)=\left\{
\begin{array}{rl}
1, \qq &\textrm{if $v$ is a terminal vertex of the edge $\be$, i.e. $t=1$ at $v$}, \\
 -1, \qq &\textrm{if $v$ is a initial vertex of the edge $\be$, i.e. $t=0$ at $v$}.
\end{array}\right.
$$

We define \emph{the
fundamental graph} $\G_F=(V_F,\cE_F)$ of the periodic graph $\Gamma$
 as a graph
on the surface $\R^d/\Z^d$ by
\[
\lb{G0} \G_F=\G/{\Z}^d\ss\R^d/\Z^d.
\]
The fundamental graph $\G_F$ has the vertex set $V_F$ and the set
$\cE_F$  of unoriented edges, which are finite.
In the space $\R^d$ we consider a coordinate system with the origin at
some point $O$ and with the basis $a_1,\ldots,a_d$. Below the coordinates of all vertices of $\G$ will be expressed  in this coordinate system.
We identify the vertices of the fundamental graph $\G_F=(V_F,\cE_F)$
with the vertices of the graph $\G=(V,\cE)$ from the set $[0,1)^d$  by
\[
\lb{V0}
V_F=[0,1)^d\cap V=\{v_1,\ldots,v_\n\},\qqq \n=\# V_F<\infty,
\]
where $\# A$ is the number of elements of the set $A$.

The metric Laplacian $\D_M$ on $L^2(\G)$ has the
decomposition into a constant fiber direct integral
\[
\lb{raz}
L^2(\G)={1\/(2\pi)^d}\int^\oplus_{\T^d}L^2(\G_F)\,d\vt ,\qqq
\mU\D_M \mU^{-1}={1\/(2\pi)^d}\int^\oplus_{\T^d}\D_M(\vt)d\vt,
\]
$\T^d=\R^d/(2\pi\Z)^d$,
for some unitary operator $\mU$.
Here the Floquet (fiber) operator $\D_M(\vt)$ acts on $y=(y_\be)_{\be\in\cE_F}\in L^2(\G_F)$
by
\[\lb{di1}
\textstyle(\D_M(\vt)y)_\be=\big(i\,{\pa\/\pa t}+\lan\t
({\bf e}),\,\vt\ran\big)^2y_{\be},\qqq (y\,''_\be)_{\be\in\cE_F}\in L^2(\G_F),
\]
see \cite{KS14c}, where
$y$ satisfies the Kirchhoff conditions \er{Dom1}; $\t
({\bf e})\in\Z^d$ is the so-called edge index, defined in subsection 4.1, $\lan\cdot\,,\cdot\ran$ denotes the standard
inner product in $\R^d$. It is more convenient for us instead of the energy parameter $E$ to introduce a new physical parameter, the momentum $z=\sqrt{E}$.
Each Floquet operator  $\D_M(\vt)$, $\vt\in\T^d$, acts on the compact graph $\G_F$ and its spectrum consists of infinitely many isolated
eigenvalues $E_n(\vt)=z_n^2(\vt)$, $n\in\N$, of finite multiplicity labeled by
\[
\label{eq.3M}
z_1^2(\vt)\leq z_2^2(\vt)\le\ldots\,.
\]
Each
$z_n^2(\cdot)$, $n\in\N$, is a real and continuous function on the
torus $\T^d$ and creates the spectral band $\s_n(\D_M)$ given by
\[
\lb{banM} \s_n(\D_M)=\big[(z_n^-)^2,(z_n^+)^2\big]=z_n^2(\T^d).
\]
Note that if
$z_n^2(\cdot)= C_n=\const$ on some set $\mB\ss\T^d$ of positive
Lebesgue measure, then the operator $\D_M$ on $\G$ has the eigenvalue
$C_n$ with infinite multiplicity. We call $C_n$ a \emph{flat band}.
The spectrum of the metric Laplace operator $\D_M$ on the periodic graph $\G$
has the form
\[
\lb{r0M}
\s(\D_M)=\bigcup_{n=1}^{\infty}\s_n(\D_M)=\s_{ac}(\D_M)\cup\s_{fb}(\D_M).
\]
Here $\s_{ac}(\D_M)$ is the absolutely continuous spectrum, which is a
union  of non-degenerated intervals from \er{banM}, and
$\s_{fb}(\D_M)$
is the set of
all flat bands (eigenvalues of infinite multiplicity). An open
interval between two neighboring non-degenerated spectral bands is
called a \emph{spectral gap}.

\subsection{Localization of spectral bands.}

Instead of the Laplacian $\D_M\ge 0$
it is convenient for us to define the momentum operator $\sqrt{\D_M}\ge 0$.
Due to Cattaneo Theorem (see Section~\ref{Sec3}) both the sets
$\s_{ac}(\sqrt{\D_M}\,)$ and $\s_{fb}(\sqrt{\D_M}\,)$
 are $2\pi$-periodic on the half-line $(0,\infty)$ and  are symmetric on the interval $(0,2\pi)$ with respect to the point $\pi$. Thus,
in order to study $\D_M$ it is sufficient to study its restriction $\Om$ on the spectral interval $[0,\pi]$ given by
\[
\lb{O1} \Om=\sqrt{\D_M}\,\chi_{[0,\pi]}(\sqrt{\D_M}\,)\;,
\]
where $\chi_A(\cdot)$ is the characteristic function of the set
$A$.
The spectrum of the operator $\Om$ on $L^2(\G)$ has  the form
\[\lb{Qr.1}
\displaystyle \s(\Om)=\bigcup_{n=1}^\n\s_n(\Om)=\s_{ac}(\Om)\cup \s_{fb}(\Om),\qqq
\s_n(\Om)=[z_n^-,z_n^+].
\]
Here $\s_{ac}(\Om)$ is a
union of non-degenerated spectral bands $\s_n(\Om)$ with $z_n^-<z_n^+\leq\pi$ and $\s_{fb}(\Om)$ is the flat band spectrum (for more details see Section \ref{Sec3}).

\

A subgraph $\G_1=(V_1,\cE_1)$ of $\G$ is called a \emph{fundamental
domain} of $\G$ if it satisfies the following conditions:

1) \emph{$\G_1=(V_1,\cE_1)$ is a finite connected graph with
  an edge set $\cE_1$ and a vertex set $V_1\supset V_F$;}

2) {\it $\G_1$ does not contain any $\Z^d$-equivalent edges;}

3) {\it $\bigcup\limits_{\mm\in\Z^d}\big(\G_1+\mm\big)=\G$.}

\no \textbf{Remark}. It is possible to remove the "convenient"
condition $V_F\ss V_1$ from the definition of the fundamental domain
and to consider a wider class of the fundamental domains. In this
case the main results still hold true, but their proof will be a bit more complicated.

The fundamental domain $\G_1$ is not uniquely defined and we fix one
of them. Let $\vk_v^1$ be the degree of the vertex $v\in V_1$ on
$\G_1$. A vertex $v\in V_1$ is called an \emph{inner} vertex of
$\G_1$, if $\vk_v=\vk_v^1$, i.e., if all its incident edges $\be\in\cE$ also belong to $\cE_1$. Denote by $V_o$ the set of all
inner vertices of $\G_1$. We define a \emph{boundary} $\pa V_1$ of
$\G_1$ by the standard identity:
\[\lb{inn}
\pa V_1=V_1\sm V_o.
\]

\no \textbf{Remark.} 1) If the graph $\G_1$ is "rather big", then the number of the inner vertices is significantly greater than the number of the boundary vertices.
If the graph $\G_1$ is "rather small", then the set $V_o$ may be empty and all vertices of $\G_1$ are the boundary vertices. But the boundary never disappears.
Some examples and discussion of the inner vertex set and the boundary see in \cite{KS14b}

2) $V_o\ss V_F$ (see Lemma 2.1 in \cite{KS14b}).

\

On the finite graph $\G_1$ we define two self-adjoint operators
$\D_{M}^1$ and $\D_{M}^o$:

1) The Neumann operator $\D_{M}^1$ on $L^2(\G_1)$ is the metric Laplacian on the graph $\G_1$, defined by \er{MLa}, \er{Dom1}.

2) The self-adjoint Dirichlet operator $\D_{M}^o$ on $f\in L^2(\G_1)$
is defined by
\[
\lb{OmD+} \D_{M}^o f=\D_{M}^1 f,\qq \textrm{where} \qq  f|_{\pa V_1}=0.
\]

Let $\Om_1$ and $\Om_o$ be the restrictions of the operators $\sqrt{\D_{M}^1}$ and $\sqrt{\D_{M}^o}$, respectively, on the spectral interval $[0,\pi]$ given by
\[
\lb{Om01} \Om_\f=(\D_M^\f)^{1/2}\chi_{[0,\pi]}\big((\D_M^\f)^{1/2}\big), \qqq \f=o,1.
\]
The spectrum of the operators $\Om_\f$, $\f=o,1$, on the finite graph $\G_1$ consists of $\n_\f$ eigenvalues, $\n_\f=\# V_\f$ and may be the additional eigenvalue $z_{\n_{\f}+1}^\f=\pi$ with an eigenfunction, vanishing at each vertex of $\G_1$.
Denote the first $\n_\f$ eigenvalues, counted
according to multiplicity, by
\[
\lb{MLDN} z_1^\f \le z_2^\f \le\ldots\le z_{\n_\f}^\f,\qqq
\n_\f=\# V_\f,\qq \f=o,1.
\]

Lled\'o and Post \cite{LP08} estimated the
position of each band $\s_n(\Om)$ for $\Om$ by
\[\lb{ELP}
\s_n(\Om)\subset J_n,\qqq n\in\N_\n,
\]
where the intervals $J_n$ have the form
\[\lb{J1M}
J_n=\ca [z_n^1,z_n^o],  & n=1,\ldots, \n_o\\[2pt]
[z_n^1,\pi],  & n=\n_o+1,\ldots,\n
\ac.
\]

The following theorem improves Lled\'o -- Post's results (see subsection 3.3).

\begin{theorem}\label{T2QG}
Each band $\s_n(\Om)$ of the operator $\Om$ acting on $L^2(\G)$
satisfies
\[\lb{ellq}
\s_n(\Om)\subset J_n\cap K_n,\qqq n=1,\ldots,\n,
\]
where the  intervals $J_n$ are defined by \er{J1M} and the intervals $K_n$ are given by
\[\lb{wJ1M}
K_n=\ca [0,z_{n+\n_1-\n}^1],  & n=1,\ldots,\n-\n_o\\[2pt]
[z_{n-\n+\n_o}^o,z_{n+\n_1-\n}^1], & n=\n-\n_o+1,\ldots,\n
\ac.
\]
\end{theorem}

\no \textbf{Remark.} 1) A graph is called \emph{bipartite} if its vertex set is
divided into two disjoint sets (called \emph{parts} of the graph)
such that each edge connects vertices from distinct parts. For a bipartite graph the interval $K_n=\zeta(J_{\n-n+1})$  for each $n\in\N_\n$, where
$\zeta(z)=\pi-z$. Thus, in this case the estimate \er{ellq} has
the form
$$
\s_n(\Om)\subset J_n\cap\zeta(J_{\n-n+1}),\qqq n\in\N_\n.
$$
Note that this result coincides with the result obtained by Lled\'o and Post in \cite{LP08}.

\begin{figure}[h]
\centering
\unitlength 1.0mm 
\linethickness{0.4pt}
\ifx\plotpoint\undefined\newsavebox{\plotpoint}\fi 

\begin{picture}(120,60)(0,0)
\put(34,28){$\scriptstyle a_2$}
\put(68.5,7.5){$\scriptstyle a_1$}
\put(34.5,8){$\scriptstyle O$}
\put(0,8){\emph{a)}}

\put(54,7.5){$\scriptstyle v_5$}
\put(54,12.5){$\scriptstyle v_1$}
\put(53.9,31.5){$\scriptstyle v_6$}
\put(54,26.5){$\scriptstyle v_2$}

\put(64.5,21.0){$\scriptstyle v_3$}
\put(41.5,21.0){$\scriptstyle v_4$}
\put(77.0,21.0){$\scriptstyle v_7$}

\put(37.5,10){\vector(1,0){35.00}}
\put(37.5,10){\vector(0,1){20.00}}
\multiput(37.5,30)(4,0){9}{\line(1,0){2}}
\multiput(72.5,10)(0,4){5}{\line(0,1){2}}
\multiput(37.5,10)(0,4){5}{\line(0,1){2}}

\put(20,10){\line(-1,1){10.00}}
\put(20,10){\line(1,1){10.00}}
\put(20,30){\line(1,-1){10.00}}
\put(20,30){\line(-1,-1){10.00}}

\put(10,20){\line(1,0){20.00}}
\put(45,20){\line(1,0){20.00}}
\put(45,20.2){\line(1,0){20.00}}
\put(45,20.1){\line(1,0){20.00}}
\put(45,19.9){\line(1,0){20.00}}
\put(80,20){\line(1,0){20.00}}

\put(10,20){\line(2,1){10.00}}
\put(30,20){\line(-2,1){10.00}}
\put(10,20){\line(2,-1){10.00}}
\put(30,20){\line(-2,-1){10.00}}
\put(20,25){\circle{1}}
\put(20,15){\circle{1}}

\put(45,20){\line(2,1){10.00}}
\put(45,20.2){\line(2,1){10.00}}
\put(45,19.9){\line(2,1){10.00}}

\put(65,20){\line(-2,1){10.00}}
\put(65,20.2){\line(-2,1){10.00}}
\put(65,19.9){\line(-2,1){10.00}}

\put(45,20){\line(2,-1){10.00}}
\put(45,20.2){\line(2,-1){10.00}}
\put(45,19.9){\line(2,-1){10.00}}

\put(65,20){\line(-2,-1){10.00}}
\put(65,20.2){\line(-2,-1){10.00}}
\put(65,19.9){\line(-2,-1){10.00}}

\put(55,25){\circle*{1.8}}
\put(55,15){\circle*{1.8}}

\put(80,20){\line(2,1){10.00}}
\put(100,20){\line(-2,1){10.00}}
\put(80,20){\line(2,-1){10.00}}
\put(100,20){\line(-2,-1){10.00}}
\put(90,25){\circle{1}}
\put(90,15){\circle{1}}

\put(10,20){\circle{1}}
\put(20,10){\circle{1}}
\put(30,20){\circle{1}}
\put(20,30){\circle{1}}

\put(55,10){\line(-1,1){10.00}}
\put(55,10.2){\line(-1,1){10.00}}
\put(55,9.8){\line(-1,1){10.00}}

\put(55,10){\line(1,1){10.00}}
\put(55,10.2){\line(1,1){10.00}}
\put(55,9.8){\line(1,1){10.00}}

\put(55,30){\line(1,-1){10.00}}
\put(55,30.2){\line(1,-1){10.00}}
\put(55,29.8){\line(1,-1){10.00}}

\put(55,30.2){\line(-1,-1){10.00}}
\put(55,30){\line(-1,-1){10.00}}
\put(55,29.8){\line(-1,-1){10.00}}

\put(30,20){\line(1,0){15.00}}

\put(45,20){\circle{1.8}}
\put(55,10){\circle{1.8}}
\put(65,20){\circle*{1.8}}
\put(55,30){\circle{1.8}}


\put(90,10){\line(-1,1){10.00}}
\put(90,10){\line(1,1){10.00}}
\put(90,30){\line(1,-1){10.00}}
\put(90,30){\line(-1,-1){10.00}}
\put(65,20.2){\line(1,0){15.00}}
\put(65,20.1){\line(1,0){15.00}}
\put(65,20){\line(1,0){15.00}}
\put(65,19.9){\line(1,0){15.00}}
\put(100,20){\line(1,0){15.00}}
\put(115,20){\circle{1}}

\put(80,20){\circle{1.8}}
\put(90,10){\circle{1}}
\put(100,20){\circle{1}}
\put(90,30){\circle{1}}
\put(20,30){\line(-1,1){10.00}}
\put(20,30){\line(1,1){10.00}}
\put(20,50){\line(1,-1){10.00}}
\put(20,50){\line(-1,-1){10.00}}

\put(10,40){\line(1,0){20.00}}
\put(45,40){\line(1,0){20.00}}
\put(80,40){\line(1,0){20.00}}

\put(10,40){\line(2,1){10.00}}
\put(30,40){\line(-2,1){10.00}}
\put(10,40){\line(2,-1){10.00}}
\put(30,40){\line(-2,-1){10.00}}
\put(20,45){\circle{1}}
\put(20,35){\circle{1}}

\put(45,40){\line(2,1){10.00}}

\put(65,40){\line(-2,1){10.00}}

\put(45,40){\line(2,-1){10.00}}

\put(65,40){\line(-2,-1){10.00}}

\put(55,45){\circle{1}}
\put(55,35){\circle{1}}

\put(80,40){\line(2,1){10.00}}
\put(100,40){\line(-2,1){10.00}}
\put(80,40){\line(2,-1){10.00}}
\put(100,40){\line(-2,-1){10.00}}
\put(90,45){\circle{1}}
\put(90,35){\circle{1}}

\put(10,40){\circle{1}}
\put(20,30){\circle{1}}
\put(30,40){\circle{1}}
\put(20,50){\circle{1}}

\put(55,30){\line(-1,1){10.00}}

\put(55,30){\line(1,1){10.00}}

\put(55,50){\line(1,-1){10.00}}

\put(55,50){\line(-1,-1){10.00}}

\put(30,40){\line(1,0){15.00}}

\put(45,40){\circle{1}}
\put(65,40){\circle{1}}
\put(55,50){\circle{1}}

\put(90,30){\line(-1,1){10.00}}
\put(90,30){\line(1,1){10.00}}
\put(90,50){\line(1,-1){10.00}}
\put(90,50){\line(-1,-1){10.00}}
\put(65,40){\line(1,0){15.00}}

\put(100,40){\line(1,0){15.00}}
\put(115,40){\circle{1}}

\put(80,40){\circle{1}}
\put(90,30){\circle{1}}
\put(100,40){\circle{1}}
\put(90,50){\circle{1}}

\end{picture}
\begin{picture}(130,40)(0,0)
\unitlength 1.0mm
\put(0,7){\vector(1,0){130.00}}
\put(0,7.1){\line(1,0){36.29}}
\put(0,6.9){\line(1,0){36.29}}
\put(0,7.2){\line(1,0){36.29}}
\put(0,6.8){\line(1,0){36.29}}
\put(36.29,6){\line(0,1){2.00}}

\put(60,7){\circle*{1}}
\put(60,9){\circle*{1}}

\put(60,7.1){\line(1,0){12.99}}
\put(60,6.9){\line(1,0){12.99}}
\put(60,7.2){\line(1,0){12.99}}
\put(60,6.8){\line(1,0){12.99}}
\put(72.99,6){\line(0,1){2.0}}

\put(77.16,7.1){\line(1,0){19.48}}
\put(77.16,6.9){\line(1,0){19.48}}
\put(77.16,7.2){\line(1,0){19.48}}
\put(77.16,6.8){\line(1,0){19.48}}
\put(77.16,6){\line(0,1){2.0}}

\put(96.64,6){\line(0,1){2.0}}

\put(0,6){\line(0,1){2.00}}
\put(120,7){\circle*{1}}
\put(-1.0,3){$\scriptstyle 0$}
\put(56.5,4){$\scriptstyle \s_2(\Om)$}
\put(56.5,1){$\scriptstyle \s_3(\Om)$}
\put(65.5,4){$\scriptstyle \s_4(\Om)$}
\put(84.5,4){$\scriptstyle \s_5(\Om)$}
\put(119.5,4){$\scriptstyle \pi$}
\put(125.0,3){$\scriptstyle \s(\Om)$}

\put(13.0,4){$\scriptstyle \s_1(\Om)$}

\put(0,18){\vector(1,0){130.00}}
\put(0,17){\line(0,1){2.00}}
\put(120,18){\circle*{1}}

\put(122.0,20){$\scriptstyle \s(\Om_o)$}
\put(60,18){\circle*{1}}
\put(43.93,18){\circle*{1}}
\put(76.01,18){\circle*{1}}

\bezier{25}(60,7)(60,15.5)(60,24)
\bezier{25}(43.93,7)(43.93,15.5)(43.93,24)
\bezier{25}(76.01,7)(76.01,15.5)(76.01,24)

\put(0,23.1){\line(1,0){60.00}}
\put(0,23.2){\line(1,0){60.00}}
\put(0,23.3){\line(1,0){60.00}}
\put(0,23.4){\line(1,0){60.00}}
\put(0,23.5){\line(1,0){60.00}}

\put(43.93,22.5){\line(1,0){16.04}}
\put(43.93,22.4){\line(1,0){16.04}}
\put(43.93,22.3){\line(1,0){16.04}}
\put(43.93,22.2){\line(1,0){16.04}}
\put(43.93,22.1){\line(1,0){16.04}}

\put(60.0,22.6){\line(1,0){15.28}}
\put(60.0,22.7){\line(1,0){15.28}}
\put(60.0,22.8){\line(1,0){15.28}}
\put(60.0,22.9){\line(1,0){15.28}}
\put(60.0,23.0){\line(1,0){15.28}}

\put(76.02,23.1){\line(1,0){20.63}}
\put(76.02,23.2){\line(1,0){20.63}}
\put(76.02,23.3){\line(1,0){20.63}}
\put(76.02,23.4){\line(1,0){20.63}}
\put(76.02,23.5){\line(1,0){20.63}}
\put(15.0,20.5){$\scriptstyle K_1=K_2$}
\put(48.0,19.5){$\scriptstyle K_3$}
\put(66.0,20.0){$\scriptstyle K_4$}
\put(85.0,20.5){$\scriptstyle K_5$}

\put(19.0,25.5){$\scriptstyle J_1$}
\put(54.0,25.5){$\scriptstyle J_2$}
\put(67.0,27.0){$\scriptstyle J_3$}
\put(99.0,26.0){$\scriptstyle J_4=J_5$}
\put(0,24.5){\line(1,0){43.93}}
\put(0,24.4){\line(1,0){43.93}}
\put(0,24.3){\line(1,0){43.93}}
\put(0,24.2){\line(1,0){43.93}}
\put(0,24.1){\line(1,0){43.93}}

\put(51.95,24.1){\line(1,0){8.02}}
\put(51.95,24.2){\line(1,0){8.02}}
\put(51.95,24.3){\line(1,0){8.02}}
\put(51.95,24.4){\line(1,0){8.02}}
\put(51.95,24.5){\line(1,0){8.02}}
\put(60,24.6){\line(1,0){60}}
\put(60,24.7){\line(1,0){60}}
\put(60,24.8){\line(1,0){60}}
\put(60,24.9){\line(1,0){60}}
\put(60,25.0){\line(1,0){60}}

\put(60,25.6){\line(1,0){16.04}}
\put(60,25.7){\line(1,0){16.04}}
\put(60,25.8){\line(1,0){16.04}}
\put(60,25.9){\line(1,0){16.04}}
\put(60,26.0){\line(1,0){16.04}}

\put(0,30){\vector(1,0){130.00}}
\put(120,30){\circle*{1}}

\put(122.0,32){$\scriptstyle \s(\Om_1)$}
\put(0,30){\circle*{1}}
\put(51.95,30){\circle*{1}}
\put(60,30){\circle*{1}}
\put(60,32){\circle*{1}}
\put(60,34){\circle*{1}}
\put(75.25,30){\circle*{1}}
\put(96.65,30){\circle*{1}}
\bezier{25}(0,30)(0,18.5)(0,7)
\bezier{25}(120,30)(120,18.5)(120,7)
\bezier{25}(51.95,30)(51.95,18.5)(51.95,7)
\bezier{10}(60,30)(60,27)(60,24)
\bezier{30}(75.25,30)(75.25,18.5)(75.25,7)
\bezier{30}(96.65,30)(96.65,18.5)(96.65,7)
\put(-9,5){\emph{b)}}

\put(0,10.5){\line(1,0){43.93}}
\put(0,10.4){\line(1,0){43.93}}
\put(0,10.3){\line(1,0){43.93}}
\put(0,10.2){\line(1,0){43.93}}
\put(0,10.1){\line(1,0){43.93}}
\put(15.0,11.5){$\scriptstyle J_1\cap K_1$}
\put(51.95,10.1){\line(1,0){8.02}}
\put(51.95,10.2){\line(1,0){8.02}}
\put(51.95,10.3){\line(1,0){8.02}}
\put(51.95,10.4){\line(1,0){8.02}}
\put(51.95,10.5){\line(1,0){8.02}}
\put(50.0,11.5){$\scriptstyle J_2\cap K_2$}
\put(60,13){\circle*{1}}
\put(56.3,14.5){$\scriptstyle J_3\cap K_3$}
\put(60,11.1){\line(1,0){15.28}}
\put(60,11.2){\line(1,0){15.28}}
\put(60,11.3){\line(1,0){15.28}}
\put(60,11.4){\line(1,0){15.28}}
\put(60,11.5){\line(1,0){15.28}}
\put(65.3,12.5){$\scriptstyle J_4\cap K_4$}
\put(76.02,10.1){\line(1,0){20.63}}
\put(76.02,10.2){\line(1,0){20.63}}
\put(76.02,10.3){\line(1,0){20.63}}
\put(76.02,10.4){\line(1,0){20.63}}
\put(76.02,10.5){\line(1,0){20.63}}
\put(82.3,11.5){$\scriptstyle J_5\cap K_5$}
\end{picture}
\caption{ \footnotesize \emph{a)} A periodic graph $\G$ and its fundamental domain $\G_1$, the vertices of $\G_1$ are big points (white and black); the edges of
$\G_1$ are marked by bold lines. The set of the inner vertices (black points) and the boundary (white points) are $V_o=\{v_1,v_2,v_3\}$ and $\partial V_1=\{v_4,v_5,v_6,v_7\}$, respectively. 
\emph{b)} Eigenvalues of the operators $\Om_1$ and $\Om_o$, the intervals $J_n$ and $K_n$, $n\in\N_5$, and their intersections, the spectrum of the operator $\Om$.} \label{ff.0.11}
\end{figure}
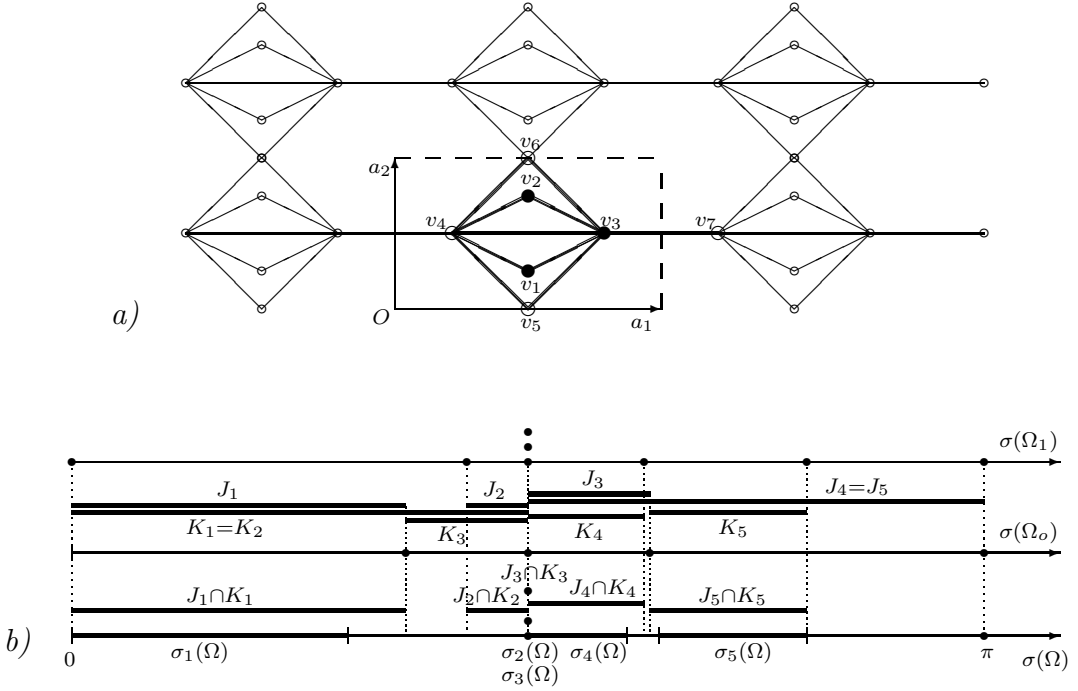

2) Theorem \ref{T2QG} estimates the positions of the spectral bands in
terms of  eigenvalues of the operators $\Om_1$ and $\Om_o$ on the fundamental domain $\G_1$. Moreover, in some cases it allows to determine the
existence of gaps and flat bands in the spectrum of the operator $\Om$. For example, for the graph shown in Fig.\ref{ff.0.11}\emph{a} the intervals $J_n\cap K_n$, $n\in\N_5$, are
shown in Fig.\ref{ff.0.11}\emph{b}. The spectrum of the operator
$\Om$ is also shown in this figure. As we can see Theorem~\ref{T2QG}
detects the flat band ${\pi\/2}$ and the
existence of all gaps in the spectrum of the operator (for more
details see Section~\ref{Sec3}).

3) Generally speaking, for distinct fundamental domains $\G_1$ the operators $\Om_1,\Om_o$, their eigenvalues and, consequently, the intervals $J_n,K_n$ are different. We number the fundamental domains $\G_1^{1},\G_1^{2},\ldots$\,.
Thus, a more precise localization of the spectral bands of the operator $\Om$ on a periodic graph $\G$ has the form
\[\lb{pl}
\s_n(\Om)\subset \bigcap_{\a}\big(J_n^{\a}\cap K_n^{\a}\big),\qqq n\in\N_\n,
\]
where $J_n^{\a},K_n^{\a}$ are the intervals, defined by \er{J1M}, \er{wJ1M}, for the fundamental domain $\G_1^{\a}$.

4) Due to the Cattaneo correspondence between the spectra of discrete and metric Laplacians, the proof of this theorem is reduced to the proof of the spectral band localization for discrete Laplacians. Moreover, we obtain this localization not only for the discrete Laplacians but also for the discrete Schr\"odinger operators with periodic potentials.

\

We present the plan of our paper. Section \ref{Sec2} is devoted to the discrete Schr\"odinger operators with periodic potentials on periodic graphs. We formulate the result about a localization of their spectral bands in terms of eigenvalues of Dirichlet and Neumann operators on a fundamental domain of the periodic
graph. In section \ref{Sec3} we prove Theorem \ref{T2QG}. In section \ref{Sec4} we prove the spectral band localization for the discrete Schr\"odinger operators on periodic graphs and estimate the Lebesgue measure of the spectrum.

\section{Localization of spectral bands for discrete Schr\"odinger operators}
\setcounter{equation}{0}\lb{Sec2}
\subsection{Discrete Schr\"odinger operators.}
Let $\ell^2(V)$ be the
Hilbert space of all square summable functions $f:V\to \C$, equipped
with the norm
$$
\|f\|^2_{\ell^2(V)}=\sum_{v\in V}|f(v)|^2<\infty.
$$
We define the self-adjoint normalized Laplacian (i.e., the Laplace operator) $\D$ on $f\in\ell^2(V)$ by
\[
\lb{DOL}
 \big(\D f\big)(v)=
 -\frac1{\sqrt{\vk_v}}\sum\limits_{(v,\,u)_e\in\cE}\frac1{\sqrt{\vk_u}}\,f(u),
 \qquad v\in V,
\]
where ${\vk}_v$ is the degree of the vertex $v\in V$ and all loops in the sum
\er{DOL} are counted twice.

We recall the basic facts (see
\cite{Ch97}, \cite{HS04}, \cite{MW89}) for both finite and periodic graphs:

{\it (i) the point $-1$ belongs to the spectrum $\s(\D)$ and  $\s(\D)$ is contained in  $[-1,1]$, i.e.,
\[\lb{mp}
-1\in\s(\D)\subset[-1,1];
\]

(ii) a graph is bipartite iff the point $1\in\s(\D)$;

(iii)  on a periodic graph the points $\pm1$ are never flat bands of $\D$.}

\

We consider the Schr\"odinger operator $H$ acting on the Hilbert
space $\ell^2(V)$ and given by
\[
\lb{Sh}
H=\D+Q,
\]
\[
\lb{Pot} \big(Q f\big)(v)=Q(v)f(v)\qqq \forall\, v\in V,
\]
where we assume that the potential $Q$ is real valued and satisfies
$$
Q(v+a_s)=Q(v), \qqq  \forall\, (v,s)\in V\ts\N_d.
$$

\subsection{The spectrum of the Schr\"odinger operator.}\label{ss23}

The discrete Schr\"odinger operator $H=\D+Q$ on $\ell^2(V)$ has  the
decomposition into a constant fiber direct integral
\[
\lb{raz}
\begin{aligned}
& \ell^2(V)={1\/(2\pi)^d}\int^\oplus_{\T^d}\ell^2(V_F)\,d\vt ,\qqq
UH U^{-1}={1\/(2\pi)^d}\int^\oplus_{\T^d}H(\vt)d\vt,
\end{aligned}
\]
$\T^d=\R^d/(2\pi\Z)^d$,     for some unitary operator $U$. Here
$\ell^2(V_F)=\C^\nu$ is the fiber space  and the Floquet $\nu\ts\nu$
matrix  $H(\vt)$ (i.e., a fiber  matrix) is given by
\[\lb{Hvt}
H(\vt)=\D(\vt)+q,\qqq q=\diag(q_1,\ldots,q_\n),\qqq \forall\,\vt\in \T^d,
\]
and $q_j$ denote the values of the potential $Q$ on the vertex set
$V_F$ by
\[\lb{pott}
Q(v_j)=q_j, \qqq j\in \N_\n=\{1,\ldots,\nu\}.
\]
The decomposition \er{raz} is standard and follows from the
Floquet-Bloch theory \cite{RS78}. The precise expression of the
Floquet matrix $\D(\vt)$ for the
Laplacian $\D$ is given by \er{l2.15'}.
Each Floquet ${\nu\ts\nu}$ matrix  $H(\vt)$, $\vt\in\T^d$, has $\n$
eigenvalues labeled by
\[
\label{eq.3}
\l_1(\vt)\leq\ldots\leq\l_{\nu}(\vt).
\]
Note that the spectrum of the Floquet matrix $H(\vt)$ does not
depend on the choice of the coordinate origin $O$.
Each
$\l_n(\cdot)$, $n\in\N_\n$, is a real and continuous function on the
torus $\T^d$ and creates the spectral band $\s_n(H)$ given by
\[
\lb{banD} \s_n(H)=[\l_n^-,\l_n^+]=\l_n(\T^d).
\]
Thus, the spectrum of the operator $H$ on the periodic graph $\G$ is
given by
\[
\lb{r0}
\s(H)=\bigcup_{\vt\in\T^d}\s\big(H(\vt)\big)=\bigcup_{n=1}^{\nu}\s_n(H).
\]
Note that if
$\l_n(\cdot)= C_n=\const$ on some set $\mB\ss\T^d$ of positive
Lebesgue measure, then the operator $H$ on $\G$ has the eigenvalue
$C_n$ with infinite multiplicity.
Thus, the spectrum of the Schr\"odinger operator $H$ on the periodic graph $\G$
has the form
\[
\lb{r0}
\s(H)=\s_{ac}(H)\cup \s_{fb}(H),
\]
where $\s_{ac}(H)$ is the absolutely continuous spectrum, which is a
union  of non-degenerated intervals from \er{banD}, and
$\s_{fb}(H)=\{\m_1,\ldots,\m_r\}$, $r<\nu$,
is the set of
all flat bands (eigenvalues of infinite multiplicity).

\subsection{Localization of spectral bands for discrete Schr\"odinger operators.}
On the fundamental domain $\G_1$ we define two self-adjoint operators
$H_1$ and $H_o$:

1) The Neumann operator $H_1$ on $\ell^2(V_1)$ is the Schr\"odinger operator on the graph $\G_1$, defined by \er{Sh}.

2) The self-adjoint Dirichlet operator $H_o$ on $f\in \ell^2(V_1)$
is defined by
\[
\lb{ShD+} H_o f=H_1 f,\qq \textrm{where} \qq  f|_{\pa V_1}=0.
\]
We will identify the Dirichlet operator $H_o$ on $f\in \ell^2(V_1)$
with the self-adjoint Dirichlet operator $H_o$ on $f\in \ell^2(V_o)$,
since $f|_{\pa V_1}=0$.

\no \textbf{Remark.} Due to the boundary conditions $f|_{\pa V_1}=0$  we call the
operator $H_o$ the Dirichlet operator.

Denote the eigenvalues of the operators $H_\f$, $\f=o,1$, counted
according to multiplicity, by
\[
\lb{LDN} \l_1^\f \le \l_2^\f \le\ldots\le \l_{\n_\f}^\f,\qqq
\n_\f=\# V_\f,\qq \f=o,1.
\]
We rewrite the sequence
$q_1,\ldots,q_\n$ defined by \er{pott} in nondecreasing order
\[
\lb{wtqn}
q^\bullet_1\le q^\bullet_2 \le\ldots \le  q^\bullet_\n.
 \]
Here
$q^\bullet_1=q_{n_1},q^\bullet_2=q_{n_2},\ldots,q^\bullet_\n=q_{n_\n}$
for some distinct numbers $n_1, {n_2},\ldots,{n_\n}\in\N_\n$.

\

Now we formulate the main result of this section about the spectral bands localization for the discrete Schr\"odinger operator.

\begin{theorem}
\label{T2} Each band $\s_n(H)$ of the operator $H=\D+Q$ acting on
$\ell^2(V)$ satisfies
\[\lb{ell}
\s_n(H)\subset \cJ_n\cap \cK_n,\qqq n=1,\ldots,\n,
\]
where the  intervals $\cJ_n, \cK_n$ are given by
\[\lb{J1}
\cJ_n=\ca [\l_n^1,\l_n^o],  & n=1,\ldots, \n_o\\[2pt]
[\l_n^1,q^\bu_n+1],  & n=\n_o+1,\ldots,\n \ac
\]
and
\[\lb{wJ1}
\cK_n=\ca [q^\bullet_n-1,\l_{n+\n_1-\n}^1],  & n=1,\ldots,\n-\n_o\\[2pt]
[\l_{n-\n+\n_o}^o,\l_{n+\n_1-\n}^1], & n=\n-\n_o+1,\ldots,\n \ac.
\]
\end{theorem}

\no \textbf{Remark.} 1) The proof of this theorem is similar to the case of the standard Schr\"odinger operator (see Theorem 1.1 in \cite{KS14b}) and differs from it only in some technical details. But for the sake of completeness we repeat it in section \ref{Sec4}.

2) Let the graph $\G$ be bipartite. If
$H=\D$, then $\cK_n=\eta(\cJ_{\n-n+1})$  for each $n\in\N_\n$, where
$\eta(z)=-z$. Thus, in this case the estimate \er{ell} has
the form
$$
\s_n(\D)\subset J_n\cap\eta(J_{\n-n+1}),\qqq n\in\N_\n.
$$

\section{\lb{Sec3} Proof of Theorem \ref{T2QG}}
\setcounter{equation}{0}

\subsection{Cattaneo Correspondence.}

Cattaneo obtained a
correspondence between the spectrum of the Laplacian $\D_M$ on  the
equilateral metric graph and the spectrum of the Laplacian $\D$ on the
corresponding discrete graph \cite{C97}.
For the sake of completeness and the reader's convenience we recall this
correspondence.

Consider the eigenvalues problem with Dirichlet boundary conditions
\[
\lb{Dp}
-y''=E y,\qqq y(0)=y(1)=0.
\]
It is known that the spectrum of this problem  is given by $\s_D=\{(\pi n)^2 : n\in\N\}$. Here $(\pi n)^2$  is the so-called  Dirichlet eigenvalue of
the problem \er{Dp}.

We formulate Cattaneo's result \cite{C97} in the form convenient for us.
 This theorem gives a basis for describing
the spectrum  of the operator $\D_M$ in terms of $\D$, and conversely.

{\bf Theorem (Cattaneo)} \emph{i) The spectrum of the operator $\sqrt{\D_M}\geq0$
on a periodic metric graph $\G$ has the form
\[
\lb{al} \s(\sqrt{\D_M}\,)=\s_{ac}(\sqrt{\D_M}\,)\cup \s_{fb}(\sqrt{\D_M}\,),
\]
\[
\lb{MAb} \s_{ac}(\sqrt{\D_M}\,)=\big\{z\in\R_+\ : \ -\cos
z\in\s_{ac}(\D)\big\},
\]
\[
\lb{MFb} \s_{fb}(\sqrt{\D_M}\,)=\big\{z\in\R_+\ : \ -\cos
z\in\s_{fb}(\D)\big\}\cup\{\pi n : n\in\N\}.
\]}

\emph{ii)  Both the sets
 $\s_{ac}(\sqrt{\D_M}\,)$ and $\s_{fb}(\sqrt{\D_M}\,)$
 are $2\pi$-periodic on the half-line $(0,\infty)$ and  are symmetric on the interval $(0,2\pi)$ with respect to the point $\pi$.}

\emph{iii) The spectrum of the operator $\Om$ on a periodic metric
graph $\G$ has  the form
\[\lb{Qr}
\begin{array}{c}
\displaystyle \s(\Om)=\bigcup_{n=1}^\n\s_n(\Om)=\s_{ac}(\Om)\cup \s_{fb}(\Om),\\[10pt]
\s_n(\Om)=[z_n^-,z_n^+], \qqq -\cos(z_n^\pm)=\l_n^\pm,\qq n\in\N_{\n}.
\end{array}
\]
Here $\s_{ac}(\D)$ is a
union of non-degenerated spectral bands $\s_n(\Om)$ with $z_n^-<z_n^+\leq\pi$. The flat band spectrum has the form
\[\lb{rel}
\s_{fb}(\Om)=\{z_1,\ldots,z_r,\pi\},
\qqq
-\cos(z_k)=\m_k\neq1,\qqq k\in\N_r.
\]
}

\no \textbf{Remark.} 1) Cattaneo considered the Laplacian $\D_M$ on connected locally finite graphs (including finite and periodic graphs). In general, some points of the Dirichlet spectrum $\s_D$ are  not eigenvalues of the Laplacian $\D_M$.

2) The relation between the spectra of $\D$ and $\sqrt{\D_M}$ is shown in
Fig.\ref{fRel}.

3) The flat bands $\pi n$,
$n\in\N$, of the operator $\sqrt{\D_M}\,$ will be called \emph{Dirichlet flat
bands}.

\setlength{\unitlength}{1.0mm}
\begin{figure}[h]
\centering
\hspace{10mm}
\unitlength 0.9mm 
\linethickness{0.4pt}
\ifx\plotpoint\undefined\newsavebox{\plotpoint}\fi 
\begin{picture}(170,80)(0,0)

\put(-5,40){\vector(1,0){170.00}}

\put(0,5){\vector(0,1){74.00}}

\put(0.2,10){\line(0,1){22.20}}
\put(-0.2,10){\line(0,1){22.20}}
\put(0.4,10){\line(0,1){22.2}}
\put(-0.4,10){\line(0,1){22.20}}

\multiput(0,32)(4,0){41}{\line(1,0){2}}
\multiput(0,45)(4,0){41}{\line(1,0){2}}
\multiput(0,53)(4,0){41}{\line(1,0){2}}
\multiput(0,66)(4,0){41}{\line(1,0){2}}

\put(0,45){\circle*{1.0}}

\put(0.4,53){\line(0,1){13.20}}
\put(0.2,53){\line(0,1){13.20}}
\put(-0.2,53){\line(0,1){13.20}}
\put(-0.4,53){\line(0,1){13.20}}

\put(-13.0,9){$\scriptstyle\l_1^-=-1$}
\put(-5.0,31){$\scriptstyle\l_1^+$}
\put(-5.0,65){$\scriptstyle\l_3^+$}
\put(-5.0,53){$\scriptstyle\l_3^-$}
\put(-4.0,37.0){$\scriptstyle 0=z_1^-$}
\put(23.0,36.8){$\scriptstyle z_1^+$}
\put(95.0,41.5){$\scriptstyle 2\pi-z_1^+$}
\put(147.0,41.5){$\scriptstyle 2\pi+z_1^+$}
\put(10.0,41.5){$\scriptstyle \s_1(\Omega)$}
\put(41.5,41.5){$\scriptstyle \s_3(\Omega)$}
\put(22.0,20.0){$\l=-\cos z$}
\put(-16.0,45){$\scriptstyle\s_2(\D)=\m_1$}
\put(-9.5,59){$\scriptstyle\s_3(\D)$}
\put(-3.5,70){$\scriptstyle 1$}
\put(-9.5,20){$\scriptstyle\s_1(\D)$}
\put(-1,10){\line(1,0){2.00}}
\put(-1,70){\line(1,0){2.00}}

\multiput(-1,10)(4,0){41}{\line(1,0){2}}
\multiput(-1,70)(4,0){41}{\line(1,0){2}}
\bezier{600}(0,10)(15,9)(31.4,40)
\bezier{600}(94.2,40)(62.8,100)(31.4,40)
\bezier{600}(94.2,40)(125.6,-20)(157,40)

\put(34.2,40.0){\circle*{1.0}}
\put(33.0,37.0){$\scriptstyle z_1$}
\put(86.0,37.){$\scriptstyle 2\pi-z_1$}
\multiput(26.9,32.5)(0,2){4}{\line(0,1){1}}
\multiput(34.2,40.0)(0,2){3}{\line(0,1){1}}
\multiput(51.2,40.7)(0,2){13}{\line(0,1){1}}
\multiput(39.2,40.1)(0,2){7}{\line(0,1){1}}

\put(62.8,40.0){\circle*{1.0}}
\put(125.6,40.0){\circle*{1.5}}
\put(62.0,37.5){$\scriptstyle \pi$}
\put(124.0,36.5){$\scriptstyle 2\pi$}

\multiput(74.4,40.7)(0,2){13}{\line(0,1){1}}
\multiput(86.4,40.1)(0,2){7}{\line(0,1){1}}
\multiput(98.7,32.5)(0,2){4}{\line(0,1){1}}
\multiput(91.4,40.0)(0,2){3}{\line(0,1){1}}
\multiput(152.4,32.5)(0,2){4}{\line(0,1){1}}
\put(74.4,40.3){\line(1,0){12}}
\put(74.4,40.2){\line(1,0){12}}
\put(74.4,40.1){\line(1,0){12}}
\put(74.4,39.9){\line(1,0){12}}
\put(74.4,39.8){\line(1,0){12}}
\put(74.4,39.7){\line(1,0){12}}
\put(91.4,40.0){\circle*{1.0}}

\put(98.7,40.3){\line(1,0){53.8}}
\put(98.7,40.2){\line(1,0){53.8}}
\put(98.7,40.1){\line(1,0){53.8}}
\put(98.7,39.9){\line(1,0){53.8}}
\put(98.7,39.8){\line(1,0){53.8}}
\put(98.7,39.7){\line(1,0){53.8}}
\put(39.2,40.3){\line(1,0){12}}
\put(39.2,40.2){\line(1,0){12}}
\put(39.2,40.1){\line(1,0){12}}
\put(39.2,39.9){\line(1,0){12}}
\put(39.2,39.8){\line(1,0){12}}
\put(39.2,39.7){\line(1,0){12}}
\put(38.0,37.0){$\scriptstyle z_3^-$}
\put(50.0,37.0){$\scriptstyle z_3^+$}

\put(0,40.3){\line(1,0){26.9}}
\put(0,40.2){\line(1,0){26.9}}
\put(0,40.1){\line(1,0){26.9}}
\put(0,39.9){\line(1,0){26.9}}
\put(0,39.8){\line(1,0){26.9}}
\put(0,39.7){\line(1,0){26.9}}

\put(162.0,36.0){$z$}

\put(-4,76.0){$\l$}

\end{picture}
\caption{\footnotesize Relation between the spectra of $\D$ and $\sqrt{\D_M}$.}
\label{fRel}
\end{figure}
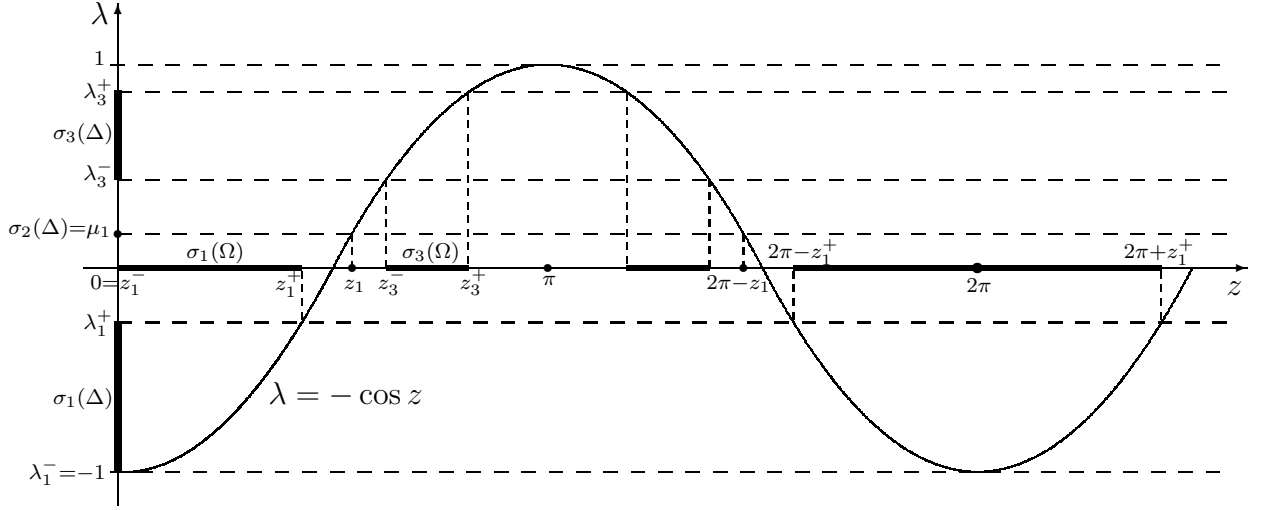

4) The number of flat bands of the operator $\Om$
is $r+1$. Flat bands $z_1,\ldots,z_r$ correspond to $r$ flat bands
of the discrete Laplacian and the flat band $\pi$ is a Dirichlet flat band.

5) For a finite graph $\G_1$ the Cattaneo correspondence between the spectra of the discrete Laplacian $H_1=\D_1$ and the operator $\Om_1$ on $\G_1$ and the similar correspondence between the spectra of these operators with Dirichlet boundary conditions on $\pa V_1$ (i.e, the spectra of the operators $H_o=\D_o$ and $\Om_o$) are given by:
\[\lb{CCF}
z\in\s(\Om_\f) \qqq \textrm{iff} \qqq -\cos(z)\in\s(H_\f), \qqq z\neq\pi, \qq \f=o,1,
\]
preserving the multiplicity of the eigenvalues (for $\f=o$ see Proposition 4.1 in \cite{LP08}). Moreover, if $1\in\s(H_1)$, i.e., the finite graph $\G_1$ is bipartite, then there exists the eigenvalue $\pi\in\s(\Om_1)$ with an eigenfunction, not vanishing at any vertex of $\G_1$ (see Lemma 4.3 in \cite{LP08}). 

\subsection{Proof of Theorem \ref{T2QG}.} Consider the spectral band $\s_n(\D)=[\l_n^-,\l_n^+]\ss[-1,1]$ of the discrete Laplacian $\D$ for some $n\in\N_\n$. Due to Cattaneo Theorem.iii the corresponding spectral band $\s_n(\Om)$ of the momentum operator $\Om$ has the form
\[\lb{f.1}
\s_n(\Om)=[z_n^-,z_n^+]\ss[0,\pi], \qq \textrm{where} \qq -\cos z_n^\pm=\l_n^\pm.
\]
Applying Theorem \ref{T2} to the spectral band $\s_n(\D)$, we obtain
\[\lb{ell'}
\s_n(\D)\ss \cJ_n\cap \cK_n,
\]
where the  intervals $\cJ_n,\cK_n\ss[-1,1]$ are given by \er{J1}, \er{wJ1} with $q^\bullet_n=0$. Since the function $\xi(\l)=\arccos(-\l)$ is an increasing bijection of the segment $[-1,1]$ onto the segment $[0,\pi]$, from \er{ell'} it follows that
\[\lb{prsb}
\xi\big(\s_n(\D)\big)\ss \xi(\cJ_n)\cap\xi(\cK_n).
\]
Due to \er{f.1}, \er{J1} and \er{wJ1}, we have
\[\lb{sOm}
\xi\big(\s_n(\D)\big)=\s_n(\Om),
\]
\[\lb{prJ1}
\xi(\cJ_n)=\ca \big[\xi(\l_n^1),\xi(\l_n^o)\big],  & n=1,\ldots, \n_o\\[2pt]
\big[\xi(\l_n^1),\pi\big],  & n=\n_o+1,\ldots,\n \ac
\]
and
\[\lb{prwJ1}
\xi(\cK_n)=\ca\big[0,\xi(\l_{n+\n_1-\n}^1)\big],  & n=1,\ldots,\n-\n_o\\[2pt]
\big[\xi(\l_{n-\n+\n_o}^o),\xi(\l_{n+\n_1-\n}^1)\big], & n=\n-\n_o+1,\ldots,\n \ac.
\]
Note that $\l_{\n_o}^o<1$, since $\s_\n\ss[\l_{\n_o}^o,\l_{\n_1}^1]\ss[-1,1]$ does not degenerate into the flat band $1$ of the discrete Laplacian $\D$.
Since the function $\xi$ is increasing and $\l_{\n_o}^o<1$, Remark 5 after Cattaneo Theorem gives
\[\lb{rel}
\xi(\l_{n}^\f)=z_{n}^\f,\qqq n\in\N_{\n_\f}, \qqq
\f=o,1.
\]
Substituting, \er{sOm} -- \er{prwJ1} into \er{prsb} and using \er{rel}, we obtain \er{ellq}. \qq \BBox

\subsection{Example 1.} Consider the operator $\Om$ on the periodic graph $\G$ shown in Fig.\ref{ff.0.11}\emph{a}. Due to Cattaneo Theorem.iii and \er{sp1} the spectrum of $\Om$ on $\G$ consists of five bands:
\[\lb{sp1q}
\textstyle\s_1(\Om)\approx[0;0{.}95],\qq \s_2(\Om)=\s_3(\Om)=\{{\pi\/2}\},\qq
\s_4(\Om)\approx[{\pi\/2};1{.}91],\qq \s_5(\Om)\approx[2{.}02;2{.}53]
\]
and the Dirichlet flat band $\pi$.

Using \er{CCF}, \er{dsp} and directly verifying that $\pi$ is the eigenvalue of the operators $\Om_1$ and $\Om_o$, we obtain the spectra of these operators
$$
\textstyle\s(\Om_1)\approx \big\{0; 1{.}36;{\pi\/2}\,;{\pi\/2}\,;{\pi\/2}\,;1{.}97;2{.}53;\pi\big\},\qq
\s(\Om_0)\approx
\big\{1{.}15;{\pi\/2}\,;1{.}99;\pi\big\}.
$$
Thus, the intervals $J_n$ and $K_n$ defined by \er{J1M}, \er{wJ1M}
and their intersections $J_n\cap K_n$,  $n\in\N_5$, have the form
$$
\begin{array}{lll}
J_1\approx[0;1{.}15], & K_1=[0,{\pi\/2}\,], \qq  & \s_1(\Om)\approx[0;0{.}95]\ss J_1\cap K_1=J_1\approx[0;1{.}15],\\[6pt]
J_2\approx[1{.}36;{\pi\/2}\,], \qq & K_2=[0,{\pi\/2}\,],  & \s_2(\Om)=\{{\pi\/2}\,\}\ss J_2\cap K_2=J_2\approx[1{.}36;{\pi\/2}\,],\\[6pt]
J_3\approx[{\pi\/2}\,;1{.}99], & K_3\approx[1{.}15;{\pi\/2}\,], & \s_3(\Om)=\{{\pi\/2}\,\}=J_3\cap K_3, \\[6pt]
J_4=[{\pi\/2}\,,\pi], & K_4\approx[{\pi\/2}\,;1{.}97],  & \s_4(\Om)\approx[{\pi\/2}\,;1{.}91]\ss J_4\cap K_4=K_4\approx[{\pi\/2}\,;1{.}97],\\[6pt]
J_5=[{\pi\/2}\,,\pi], & K_5\approx[1{.}99;2{.}53], \qq  & \s_5(\Om)\approx[2{.}02;2{.}53]\approx J_5\cap K_5=K_5\approx[1{.}99;2{.}53].
\end{array}
$$

Theorem \ref{T2QG} determines the existence of two gaps and the flat band ${\pi\/2}$
(see Fig.\ref{ff.0.11}\emph{b}). The intersection of the intervals $J_n$ and $K_n$, $n=3,4,5$, gives a more
precise estimate of the spectral band $\s_n(\Om)$ than one interval
$J_n$. Moreover, for $n=4,5$
the estimate $\s_n(\Om)\subset J_n$ gives the upper bound
$z_n^+\leq\pi$ that is trivial. But using \er{ellq} we obtain more accurate
estimates for the spectral bands. Note that the third spectral band (that degenerates into the flat band) is detected precisely.

\no \textbf{Remark.} The Lebesgue measure $|\s(\Om)|$
and $|\s(\D)|$ of the spectrum of the operator $\Om$ and the discrete Laplacian $\D$, respectively, satisfies
\[
\lb{eq.7Q}
\textstyle |\s(\Om)|\leq{\pi\/\sqrt{2}}\,|\s(\D)|^{1\/2},
\]
see Theorem 1.1.ii in \cite{KS14a}. For the graph shown in Fig.\ref{ff.0.11}\emph{a}, using the estimates \er{eq.7Q} and \er{est33}, we obtain
\[
\textstyle |\s(\Om)|\leq{\pi\/\sqrt{2}}\,|\s(\D)|^{1/2}\leq
{\pi\/\sqrt{2}}\,
\Big(\sum\limits_{n=1}^{5}|\s_n(\D)|\Big)^{1/2}\approx
2{.}81.
\]
Finally, we note that \er{sp1q} yields
$$
\textstyle\sum\limits_{n=1}^5|\s_n(\Om)|\approx(0{.}95-0)+
(1{.}91-{\pi\/2})+(2{.}53-2{.}02)\approx1{.}80.
$$

\section{\lb{Sec4} Results for discrete Schr\"odinger operators}
\setcounter{equation}{0}

\subsection{The Floquet matrix for the discrete Schr\"odinger operator.} We need to introduce the two oriented edges $(u,v)$ and $(v,u)$ for
each unoriented edge $(u,v)_e\in \cE$: the oriented edge starting at
$u\in V$ and ending at $v\in V$ will be denoted as the ordered pair
$(u,v)$. We denote the sets of all oriented edges of the graph $\G$
and the fundamental graph $\G_F$ by $\cA$ and $\cA_F$, respectively.

We introduce {\it an edge index}, which is important to study the
spectrum of Schr\"odinger operators on periodic graphs. For any
$v\in V$ the following unique representation holds true:
\[
\lb{Dv} v=[v]+\tilde v, \qquad [v]\in\Z^d,\qquad \tilde v\in
V_F\subset[0,1)^d.
\]
In other words, each vertex $v$ can be represented uniquely as the
sum  of an integer part $[v]\in \Z^d$ and a fractional part $\tilde
v$ that is a vertex of $V_F$ defined in \er{V0}. For any oriented
edge $\be=(u,v)\in\cA$ we define {\bf the edge "index"}  $\t({\bf
e})$ as the integer vector
\[
\lb{in} \t({\bf e})=[v]-[u]\in\Z^d,
\]
where due to \er{Dv} we have
$$
u=[u]+\tilde{u},\qquad v=[v]+\tilde{v}, \qquad [u],
[v]\in\Z^d,\qquad \tilde{u},\tilde{v}\in V_F.
$$

If $\be=(u,v)$ is an oriented edge of the graph $\G$, then by the
definition of the fundamental graph there is an oriented edge
$\tilde\be=(\tilde u,\tilde v\,)$ on $\G_F$. For the edge
$\tilde\be\in\cA_F$ we define the edge index $\t(\tilde{\bf e})$ by
\[
\lb{inf} \t(\tilde{\bf e})=\t(\be).
\]
In other words, edge indices of the fundamental graph $\G_F$  are
induced by edge indices of the periodic graph $\G$. The edge
indices, generally speaking, depend on the choice of the coordinate
origin $O$ and the periods $a_1,\ldots,a_d$ of the graph $\G$. But in a fixed coordinate system the index of the
fundamental graph edge is uniquely determined by \er{inf}, since
$$
\t(\be+\mm)=\t(\be),\qqq \forall\, (\be,\mm)\in\cA \ts \Z^d.
$$

The Schr\"odinger operator $H=\D+Q$ acting on $\ell^2(V)$ has the decomposition
into a constant fiber direct integral \er{raz},
where the Floquet $\nu\ts\nu$
matrix  $H(\vt)$ has the form \er{Hvt}.
The Floquet matrix
$\D(\vt)=\{\D_{jk}(\vt)\}_{j,k=1}^\n$ for the Laplacian $\D$ is given by
\medskip
\[
\label{l2.15'}
\D_{jk}(\vt )=\ca {-1\/\sqrt{\vk_j\vk_k}}
\sum\limits_{{\bf e}=(v_j,\,v_k)\in{\cA}_F}e^{\,i\lan\t
({\bf e}),\,\vt\ran }, \qq &  {\rm if}\  \ (v_j,v_k)\in \cA_F \\
\qqq 0, &  {\rm if}\  \ (v_j,v_k)\notin \cA_F
\ac,
\]
see \cite{KS13}, where $\vk_j$ is the degree of $v_j$ and $\lan\cdot\,,\cdot\ran$ denotes the standard
inner product in $\R^d$. This explicit expression for the Floquet matrix is very important to prove our main results.

\subsection{Proof of Theorem \ref{T2}.} We need the following simple fact (see Theorem 4.3.1 in
\cite{HJ85}). \emph{Let $A,B$ be self-adjoint $\nu\ts\nu$ matrices.
Denote by $\l_1(A)\leq\ldots\leq\l_\n(A)$,
$\l_1(B)\leq\ldots\leq\l_\n(B)$ the eigenvalues of  $A$ and $B$,
respectively, arranged in increasing order, counting multiplicities.
Then we have}
\[
\lb{MP} \l_n(A)+\l_1(B)\leq\l_n(A+B)\leq\l_n(A)+\l_\n(B), \qqq
\forall\, n\in\N_\n.
\]

Inequalities \er{MP} and the basic fact \er{mp} give that the eigenvalues of the Floquet matrix $H(\vt)$ for the
Schr\"odinger operator $H=\D+Q$, satisfy
\[
\lb{qq}
q_n^\bullet-1\le\l_n(\vt)\le q_n^\bullet+1,\qqq \forall\, (\vt,n)\in\T^d\ts\N_\n.
\]

Since $V_o\ss V_F$, without loss of generality we may assume that the set $V_o$ of the inner vertices of the graph $\G_1=(V_1,\cE_1)$ has the form
$$
V_o=\{v_1,\ldots,v_{\n_o}\}.
$$
We denote the equivalence classes from $V_1/\Z^d$ by
\[\lb{Vj}
\cZ_j\equiv\cZ(v_j)=\big(\{v_j\}+\Z^d\big)\cap V_1, \qqq j\in\N_\n.
\]
Note that $\cZ_j=\{v_j\}$ for all $j\in\N_{\n_o}$.

The Neumann operator $H_1$ on the graph $\G_1$ is equivalent
to the $\n_1\ts\n_1$ self-adjoint matrix
$H_1=\{H^1_{jk}\}_{j,k=1}^{\n_1}$ given by
\[
\lb{H+} H_1=\D_1+q^1,\qqq q^1=\diag(q_1^1,\ldots,q_{\n_1}^1),
\]
where $q^1_k=q_j$, if $v_k\in\cZ_j$, $k\in\N_{\n_1}$, $j\in\N_\n$,
and the matrix
$\D_1=\{\D_{jk}^1\}_{j,k=1}^{\n_1}$ has the form
\[
\label{l2.16}
\textstyle\D^1_{jk}=-{\vk^1_{jk}\/\sqrt{\vk^1_j\vk^1_k}}\,.
\]
Here $\vk^1_j$ is the degree of the vertex $v_j\in V_1$ on the graph
$\G_1$; $\vk^1_{jk}\ge 1$ is the multiplicity of the edge
$(v_j,v_k)\in\cE_1$ and  $\vk^1_{jk}=0$ if $(v_j,v_k)\notin\cE_1$.

The Dirichlet operator $H_o$ is described by the $\n_o\ts\n_o$
self-adjoint matrix $H_o=\{H_{jk}^o\}_{j,k=1}^{\n_o}$ with entries
\[\lb{H-}
H_{jk}^o=H_{jk}^1\qqq
\textrm{for all}\qqq j,k\in\N_{\n_o}.
\]

Recall that
\[\lb{rc}
\vk_j^1=\vk_j  \qq \textrm{ for all } \qq j\in\N_{\n_o}.
\]

We recall well-known facts.

\emph{Denote by $\l_1(A)\leq\ldots\leq\l_\n(A)$ the eigenvalues of a
self-adjoint $\n\ts\n$ matrix $A$, arranged in increasing order,
counting multiplicities. Each $\l_n$ satisfies the minimax principle:}
\[
\lb{CF1} \l_n(A)=\min_{S_n\subset\C^\n}\max_{\|x\|=1 \atop x\in
S_n}\lan Ax,x\ran,
\]
\[
\lb{CF2} \l_n(A)=\max_{S_{\n-n+1}\subset\C^\n}\min_{\|x\|=1 \atop
x\in S_{\n-n+1}}\lan Ax,x\ran,
\]
\emph{where $S_n$ denotes a subspace of dimension $n$ and the outer
optimization is over all subspaces of the indicated dimension} (see
p.180 in \cite{HJ85}).

First, for each $\vt\in\T^d$ we define the $\n$-dimensional subspace $Y_\vt\ss\C^{\n_1}$ by
\[
\begin{aligned}
\lb{Yt}
Y_\vt=\big\{x=(x_k)_{k=1}^{\n_1}\in\C^{\n_1} : \forall \,
k=\n+1,\ldots,\n_1 \qq x_k=\textstyle\sqrt{{\vk_k^1\/\vk_j^1}}\, e^{\,i\lan
v_k-v_j,\,\vt\ran}\,x_j\,, \\
\textrm{where $j=j(k)\in\N_\n$ is such that $v_k\in\cZ_j$}\big\}.
\end{aligned}
\]
Note that $j=j(k)$ in \er{Yt} is uniquely defined for each $k=\n+1,\ldots,\n_1$.
Let $1\leq n\leq\n$. Using \er{CF1} and \er{CF2} we write
\[
\lb{CF3'} \l_j^1=\min_{S_j\ss \C^{\n_1}} \max_{\|x\|=1 \atop
{x\in S_j}}\lan H_1x,x\ran\geq\min_{S_j\ss \C^{\n_1}} \max_{\|x\|=1
\atop {x\in S_j\cap Y_\vt}}\lan H_1x,x\ran, \qq j=n+\n_1-\n,
\]
\[
\lb{CF3} \l_n^1=\max_{S_k\ss\C^{\n_1}}\min_{\|x\|=1 \atop {x\in
S_k}}\lan H_1x,x\ran\leq\max_{S_k\ss\C^{\n_1}}\min_{\|x\|=1 \atop
x\in S_k\cap Y_\vt}\lan H_1x,x\ran, \qq k=\n_1-n+1,
\]
where $S_j$ denotes a subspace of dimension $j$. For $x\in Y_\vt$ we have
\[\lb{CF5.11}
\lan H_1x,x\ran=\sum_{j,k=1}^{\n_1}H^1_{jk}\,\bar
x_j\,x_k=\sum_{j=1}^{\n_1}q_j^1|x_j|^2-
\sum_{j,k=1}^{\n_1}{\vk_{jk}^1\/\sqrt{\vk_j^1\vk_k^1}}\,\bar
x_j\,x_k,
\]
where
\[
\lb{CF5.22}
\sum_{j=1}^{\n_1}q^1_j|x_j|^2=
\sum_{j=1}^{\n_o}q_j|x_j|^2+\sum_{j=\n_o+1}^{\n}
{q_j\sum\limits_{v\in
\cZ_j}\vk_v^1\/\vk_j^1}\,|x_j|^2=\sum_{j=1}^{\n_o}q_j|x_j|^2
+\sum_{j=\n_o+1}^{\n}q_j\, {\vk_j\/\vk_j^1}\,|x_j|^2,
\]
and
\[\lb{CF5.33}
\sum_{j,k=1}^{\n_1}{\vk_{jk}^1\/\sqrt{\vk_j^1\vk_k^1}}\,\bar
x_j\,x_k=
\sum_{j,k=1}^{\n}{1\/\sqrt{\vk_j^1\vk_k^1}}\sum\limits_{{\bf
e}=(v_j,\,v_k)\in\cA_F}e^{\,i\lan\t ({\bf e}),\,\vt\ran }\,\bar
x_j\,x_k.
\]
In \er{CF5.22} we have used the identity
\[\lb{deg}
\sum\limits_{v\in\cZ_j}\vk_v^1=\vk_j.
\]
We introduce the new vector
\[\lb{vy}
y=(y_j)_{j=1}^\n\in \C^\n,\qqq y_j=x_j\sqrt{{\vk_j\/\vk_j^1}}\,, \qqq j\in\N_\n.
\]
Since $\vk_j^1=\vk_j $ for $1\leq j\leq \n_o$, we have $y_j=x_j$, $j\in\N_{\n_o}$, and, using \er{deg}, for
$x\in Y_\vt$ we have
\[
\begin{aligned}\lb{CF6}
\|x\|^2=\sum_{j=1}^{\n_1}|x_j|^2=\sum_{j=1}^{\n_o}|x_j|^2+\sum_{j=\n_o+1}^{\n}
{\sum\limits_{v\in
\cZ_j}\vk_v^1\/\vk_j^1}|x_j|^2\\=\sum_{j=1}^{\n_o}|x_j|^2+\sum_{j=\n_o+1}^{\n}
{\vk_j\/\vk_j^1}\,|x_j|^2=\sum_{j=1}^{\n}|y_j|^2=\|y\|^2.
\end{aligned}
\]
Combining \er{CF5.11} -- \er{CF5.33} for $x\in Y_\vt$, \er{vy} and the
definition of  $H(\vt)$ in \er{Hvt}, we obtain
\[\lb{CF5.111}
\lan H_1x,x\ran=\sum_{j=1}^{\n}q_j|y_j|^2
-\sum_{j,k=1}^{\n}{1\/\sqrt{\vk_j\vk_k}}\sum\limits_{{\bf
e}=(v_j,\,v_k)\in\cA_F}e^{\,i\lan\t ({\bf e}),\,\vt\ran }\,\bar
y_j\,y_k=\lan H(\vt)y,y\ran.
\]
This, \er{CF3'}, \er{CF3}, \er{CF6} and the minimax principle \er{CF1},
\er{CF2} yield for  $1\leq n\leq\n$:
\[
\lb{CF10'} \l_{n+\n_1-\n}^1\geq\min_{S_{n}\ss \C^{\n}}
\max_{\|y\|=1 \atop {y\in S_{n}}}\lan H(\vt)y,y\ran=\l_n(\vt)\,,
\]
\[
\lb{CF10} \l_n^1\leq\max_{S_{\n-n+1}\ss\C^{\n}} \min_{\|y\|=1
\atop {y\in S_{\n-n+1}}}\lan H(\vt)y,y\ran=\l_n(\vt)\,.
\]

Second, let $X=\{x\in\C^{\n} : x_{\n_o+1}=\ldots=x_\n=0\}$ be
the $\n_o$-dimensional subspace of $\C^{\n}$ and
let $1\leq n\leq\n_o$. Using \er{CF1} and \er{CF2} we write
\[
\lb{CF3.1'} \l_j(\vt)=\min_{S_j\ss \C^{\n}} \max_{\|x\|=1 \atop
{x\in S_j}}\lan H(\vt)x,x\ran\geq\min_{S_j\ss \C^{\n}} \max_{\|x\|=1
\atop {x\in S_j\cap X}}\lan H(\vt)x,x\ran,\qq
j=n+\n-\n_o\,,
\]
\[
\lb{CF3.1} \l_n(\vt)=\max_{S_k\ss \C^{\n}} \min_{\|x\|=1 \atop {x\in
S_k}}\lan H(\vt)x,x\ran\leq \max_{S_k\ss \C^{\n}} \min_{\|x\|=1
\atop {x\in S_k\cap X}}\lan H(\vt)x,x\ran, \qq k=\n-n+1.
\]
For $x\in X$ we have
\[\lb{CF5.1}
\lan H(\vt)x,x\ran=\sum_{j,k=1}^{\n}H_{jk}(\vt)\,\bar
x_j\,x_k=\sum_{j,k=1}^{\n_o}H_{jk}^o,\bar x_j\,x_k=\lan H_ox,x\ran,
\]
\[\lb{CF6.1}
\|x\|=\sum_{j=1}^{\n}|x_j|^2=\sum_{j=1}^{\n_o}|x_j|^2.
\]
Then for $1\leq n\leq\n_o$ we may rewrite the inequalities \er{CF3.1'}, \er{CF3.1} in the
form
\[
\lb{CF7.1'} \l_{n+\n-\n_o}(\vt)\geq\min_{S_{n}\ss \C^{\n_o}}
\max_{\|x\|=1\atop {x\in S_{n}}}\lan H_ox,x\ran=\l_n^o,
\]
\[
\lb{CF7.1} \l_n(\vt)\leq\max_{S_{\n_o-n+1}\ss \C^{\n_o}}
\min_{\|x\|=1 \atop {x\in S_{\n_o-n+1}}}\lan H_ox,x\ran=\l_n^o.
\]

Combining \er{CF10} and \er{CF7.1} and using \er{qq}, we obtain for all $\vt\in\T^d$:
\[\lb{FP}
\begin{array}{ll}
\l_n(\vt)\in [\l_n^1,\l_n^o]=\cJ_n,\qqq & n=1,\ldots,\n_o\,, \\[6pt]
\l_n(\vt)\in  [\l_n^1, q_n^\bullet+1]=\cJ_n,\qqq & n=\n_o+1,\ldots,\n.
\end{array}
\]

Similarly, from \er{CF10'} and \er{CF7.1'} we obtain
\[\lb{FP1}
\begin{array}{ll}
\l_n(\vt)\in [q_n^\bu-1, \l_{n+\n_1-\n}^1]=\cK_n,\qqq & n=1,\ldots,\n-\n_o\,, \\[6pt]
\l_n(\vt)\in [\l_{n+\n_o-\n}^o, \l_{n+\n_1-\n}^1]=\cK_n,\qqq &
n=\n-\n_o+1,\ldots,\n,
\end{array}
\]
for all $\vt\in\T^d$. The relations \er{FP} and
\er{FP1} prove \er{ell}. \qq \BBox

\

Now we estimate the total length of all spectral bands of
$H$.

\begin{theorem}
\label{TT2D} The total length of all spectral bands $\s_n(H)$, $n\in\N_\n$, of $H$
satisfies
\[\lb{est2}
\sum_{n=1}^{\n}|\s_n(H)|
\le\sum_{n=1}^{\n-\n_o}\big(\l_{\n_1-(\n-\n_o)+n}^1-\l_n^1\big),\qq if \qq
\n\geq2\n_o,\;
\]
\[\lb{est3}
\sum_{n=1}^{\n}|\s_n(H)|
\le\sum_{n=1}^{\n-\n_o}\big(\l_n^o-\l_{(2\n_o-\n)+n}^o\big)-
\sum_{n=1}^{\n_o}\big(\l_n^1-\l_{(\n_1-\n_o)+n}^1\big)\qqq if \qqq
\n<2\n_o.
\]
\end{theorem}

\no {\bf Proof.}  Let $\n>2\n_o$.
Then, using \er{J1} and \er{wJ1}, we have
\begin{multline}\lb{bol}
\sum_{n=1}^{\n}|\s_n(H)|\le
\sum_{n=1}^{\n_o}(\l_n^o-\l_n^1)+\sum_{n=\n_o+1}^{\n-\n_o}(\l_{n+\n_1-\n}^1-
\l_n^1)+\sum_{n=\n-\n_o+1}^{\n}(\l_{n+\n_1-\n}^1-
\l_{n-\n+\n_o}^o)\\=
\sum_{n=\n_o+1}^{\n}\l_{n+\n_1-\n}^1-
\sum_{n=1}^{\n-\n_o}\l_{n}^1=\sum_{n=1}^{\n-\n_o}\big(\l_{\n_1-(\n-\n_o)+n}^1-
\l_n^1\big).
\end{multline}
Similarly, if $\n=2\n_o$, then the formulas \er{J1} and \er{wJ1} give
\[\lb{equ}
\sum_{n=1}^{\n}|\s_n(H)|\le
\sum_{n=1}^{\n_o}(\l_n^o-\l_n^1)+\sum_{n=\n_o+1}^{\n}(\l_{n+\n_1-\n}^1-
\l_{n-\n+\n_o}^o)=\sum_{n=1}^{\n_o}\big(\l_{\n_1-(\n-\n_o)+n}^1-
\l_n^1\big).
\]
The estimates \er{bol} and \er{equ} give \er{est2}.

Now let $\n<2\n_o$. Then, using \er{J1} and \er{wJ1}, we have
\begin{multline*}
\sum_{n=1}^{\n}|\s_n(H)|\le
\sum_{n=1}^{\n-\n_o}(\l_n^o-\l_n^1)+\sum_{n=\n-\n_o+1}^{\n_o}(\l_{n+\n_1-\n}^1-
\l_n^1)+\sum_{n=\n_o+1}^{\n}(\l_{n+\n_1-\n}^1-
\l_{n-\n+\n_o}^o)\\=\sum_{n=1}^{\n-\n_o}\l_n^o-\sum_{n=\n_o+1}^{\n}
\l_{n-\n+\n_o}^o-\sum_{n=1}^{\n_o}\l_n^1+\sum_{n=\n-\n_o+1}^{\n}\l_{n+\n_1-\n}^1\\
=\sum_{n=1}^{\n-\n_o}\big(\l_n^o-\l_{(2\n_o-\n)+n}^o\big)-
\sum_{n=1}^{\n_o}\big(\l_n^1-\l_{(\n_1-\n_o)+n}^1\big).
\end{multline*}
Thus, the estimate \er{est3} has also been proved.
\qq \BBox

\subsection{Example 2.} Consider the Laplacian $H=\D$ on the periodic graph $\G$ shown in Fig.\ref{ff.0.11}\emph{a}. The set of the fundamental graph vertices is $V_F= \{v_1,v_2,v_3,v_4,v_5\}$. For each $\vt=(\vt_1,\vt_2)\in\T^2$ the matrix $\D(\vt)$ defined by (\ref{l2.15'})
has the form
\begin{equation}\label{z2}
\D(\vt)=\left(
\begin{array}{ccccc}
   0 & 0 & {-1\/\sqrt{12}} & {-1\/\sqrt{12}} & 0\\ [6pt]
   0 & 0 & {-1\/\sqrt{12}} & {-1\/\sqrt{12}} & 0\\ [6pt]
   {-1\/\sqrt{12}} & {-1\/\sqrt{12}} & 0 & -{1+e^{i\vt_1}\/6} & -{1+e^{i\vt_2}\/\sqrt{24}} \\ [6pt]
   {-1\/\sqrt{12}} & {-1\/\sqrt{12}} & -{1+e^{-i\vt_1}\/6} & 0 & -{1+e^{i\vt_2}\/\sqrt{24}} \\ [6pt]
   0 & 0 & -{1+e^{-i\vt_2}\/\sqrt{24}} & -{1+e^{-i\vt_2}\/\sqrt{24}} & 0
\end{array}\right).
\end{equation}
The characteristic polynomial of $\D(\vt)$ is given by
$$
\textstyle\det(\D(\vt)-\l\1_5)=\l^2\Big(-\l^3+\big({1\/18}c_1+
{1\/6}c_2+{5\/9}\big)\l-{1\/12}c_1
-{1\/36}c_1c_2-{1\/36}c_2-{1\/12}\,\Big),
$$
$\1_5$ is the $5\ts5$ identity matrix, $c_1=\cos\vt_1$, $c_2=\cos\vt_2$. The spectrum of the Laplacian $\D$ on the periodic graph $\G$
consists of five bands:
\[\lb{sp1}
\textstyle\s_1(\D)\approx[-1;-0{.}58],\qq \s_2(\D)=\s_3(\D)=\{0\},\qq
\s_4(\D)\approx[0;0{.}33],\qq \s_5(\D)\approx[0{.}43;0{.}82].
\]

The fundamental domain $\G_1=(V_1,\cE_1)$  shown in
Fig.\ref{ff.0.11}\emph{a}
has the vertex set $V_1$ given by
$$
V_1=\{v_1,v_2,v_3,v_4,v_5,v_6=v_5+a_2,v_7=v_4+a_1\}.
$$
The set of the inner vertices $V_o$ and the boundary
$\partial V_1$ of $\G_1$ have the form
$$
V_o=\{v_1,v_2,v_3\},\qqq \pa V_1=\{v_4,v_5,v_6,v_7\}.
$$
The matrices $H_1$ and $H_o$, defined by \er{H+} -- \er{H-}, in this
case have the form
$$
H_1=\left(
\begin{array}{ccccccc}
0 & 0  & {-1\/\sqrt{12}} & {-1\/\sqrt{10}} & 0 & 0 & 0\\[6pt]
0 & 0  & {-1\/\sqrt{12}} & {-1\/\sqrt{10}} & 0 & 0 & 0 \\[6pt]
{-1\/\sqrt{12}} & {-1\/\sqrt{12}}  & 0 & {-1\/\sqrt{30}} & {-1\/\sqrt{12}} & {-1\/\sqrt{12}} & {-1\/\sqrt{6}}\\[6pt]
{-1\/\sqrt{10}} & {-1\/\sqrt{10}} & {-1\/\sqrt{30}} & 0 & {-1\/\sqrt{10}} & {-1\/\sqrt{10}} & 0 \\[6pt]
0 & 0 & {-1\/\sqrt{12}} &  {-1\/\sqrt{10}} & 0 & 0 & 0 \\[6pt]
0 & 0 & {-1\/\sqrt{12}} &  {-1\/\sqrt{10}} & 0 & 0 & 0 \\[6pt]
0 & 0 & {-1\/\sqrt{6}} &  0 & 0 & 0 & 0 \\[6pt]
\end{array}\right),\qqq H_o=\left(
\begin{array}{ccc}
0 & 0  & {-1\/\sqrt{12}} \\[6pt]
0 & 0  & {-1\/\sqrt{12}}  \\[6pt]
{-1\/\sqrt{12}}  & {-1\/\sqrt{12}}  & 0 \\[6pt]
\end{array}\right).
$$
The spectra of the operators $H_1$ and $H_o$ are
\[\lb{dsp}
\textstyle\s(H_1)\approx \big\{-1; -0{.}21;0;0;0;0{.}39;0{.}82\big\},\qq
\s(H_0)= \big\{-{1\/\sqrt{6}}\,;0\,;{1\/\sqrt{6}}\,\big\}\approx
\big\{-0{.}41;0;0{.}41\big\}.
\]
Thus, the intervals $\cJ_n$ and $\cK_n$ defined by \er{J1}, \er{wJ1}
and their intersections $\cJ_n\cap\cK_n$,  $n\in\N_5$, have the form
$$
\begin{array}{lll}
\cJ_1\approx[-1;-0{.}41], & \cK_1=[-1,0], \qq  & \s_1(\D)\approx[-1;-0{.}58]\ss \cJ_1\cap\cK_1=\cJ_1\approx[-1;-0{.}41],\\[6pt]
\cJ_2\approx[-0{.}21;0], \qq & \cK_2=[-1,0],  & \s_2(\D)=\{0\}\ss \cJ_2\cap\cK_2=\cJ_2\approx[-0{.}21;0],\\[6pt]
\cJ_3\approx[0;0{.}41], & \cK_3\approx[-0{.}41;0], & \s_3(\D)=\{0\}=\cJ_3\cap\cK_3, \\[6pt]
\cJ_4=[0,1], & \cK_4\approx[0;0{.}39],  & \s_4(\D)\approx[0;0{.}33]\ss \cJ_4\cap\cK_4=\cK_4\approx[0;0{.}39],\\[6pt]
\cJ_5=[0,1], & \cK_5\approx[0{.}41;0{.}82], \qq  & \s_5(\D)\approx[0{.}43;0{.}82]\approx \cJ_5\cap\cK_5=\cK_5\approx[0{.}41;0{.}82].
\end{array}
$$

\no \textbf{Remark.} For the graph shown in Fig.\ref{ff.0.11}\emph{a} $6=2\n_o>\n=5$ and the estimate \er{est3} has the form
\[
\begin{aligned}\label{est33}
\sum_{n=1}^{5}|\s_n(H)|
\le\sum_{n=1}^{2}\big(\l_n^o-\l_{n+1}^o\big)-
\sum_{n=1}^{3}\big(\l_n^1-\l_{n+4}^1\big)\\=(\l_1^o-\l_3^o)-(\l_1^1+\l_2^1+\l_3^1)+
(\l_5^1+\l_6^1+\l_7^1)
\approx-0{.}82+1{.}21+1{.}21=1{.}60.
\end{aligned}
\]
Finally, we note that  \er{sp1} yields
$$
\sum_{n=1}^5|\s_n(H)|\approx(-0{.}58+1)+(0{.}33-0)+(0{.}82-0{.}43)=1{.}14.
$$

\medskip

\no\textbf{Acknowledgments.}
\footnotesize
Various
parts of this paper were written during Evgeny Korotyaev's stay  in
the Mathematical Institute of Tsukuba University, Japan. He is grateful to
the institutes for the hospitality. His study was partly supported by the RFFI grant
No. 11-01-00458 and by the project SPbSU No. 11.38.215.2014.


\begin{thebibliography}{9999}
\setlength{\itemsep}{-\parskip}
\footnotesize

\bibitem[C97]{C97} Cattaneo, C. The spectrum of the continuous Laplacian
 on a graph, Monatsh. Math. 124 (1997), 215--235.

\bibitem[Ch97]{Ch97} Chung, F. Spectral graph theory, AMS,
 Providence, RI, 1997.


\bibitem[HS04]{HS04} Higuchi, Y. Shirai, T. Some spectral and
 geometric properties for infinite graphs, AMS Contemp. Math. 347 (2004), 29--56.

\bibitem[HJ85]{HJ85} Horn, R; Johnson, C. Matrix analysis.
Cambridge University Press, 1985.

\bibitem[KS13]{KS13} Korotyaev, E.; Saburova, N. Spectral estimates for
normalized Laplacian and its perturbations on periodic discrete
graphs, preprint 2013.

\bibitem[KS14a]{KS14a} Korotyaev, E.; Saburova, N. Estimates of bands for
Laplacians on periodic equilateral metric
graphs, preprint 2014.

\bibitem[KS14b]{KS14b} Korotyaev, E.; Saburova, N. Spectral band localization  for Schr\"odinger operators on discrete periodic graphs, preprint 2014.

\bibitem[KS14c]{KS14c} Korotyaev, E.; Saburova, N. Laplacians on periodic metric graphs, in preparation.

\bibitem[LP08]{LP08} Lled\'o, F.; Post, O. Eigenvalue bracketing for
 discrete and metric graphs, J. Math. Anal. Appl. 348 (2008), 806--833.

\bibitem[MW89]{MW89} Mohar, B.; Woess, W. A survey on spectra of infinite graphs,
Bull. London Math. Soc., 21 (1989), 209--234.

\bibitem[RS78]{RS78} Reed, M.; Simon, B. Methods of modern mathematical
physics, vol.IV, Analysis of operators, Academic Press, New York, 1978.

\end{thebibliography}
\end{document}